\documentclass[12pt]{article}
\usepackage{amssymb,amsfonts,amsmath, psfrag,eepic,colordvi,graphicx,epsfig,ytableau}
\usepackage{amssymb,latexsym,graphics,array}
\usepackage{MnSymbol}
\usepackage[enableskew]{youngtab}
\usepackage{float}
\usepackage{tikz}
\topskip  
\newcommand{\rmnum}[1]{\romannumeral #1}

\numberwithin{equation}{section}
\newtheorem{theorem}{Theorem}[section]

\newtheorem{conjecture}[theorem]{Conjecture}

\newtheorem{lemma}[theorem]{Lemma}

\newtheorem{example}[theorem]{Example}

\begin{document}
	\parskip 6pt
	
	\pagenumbering{arabic}
	\def\sof{\hfill\rule{2mm}{2mm}}
	\def\ls{\leq}
	\def\gs{\geq}
	\def\SS{\mathcal S}
	\def\qq{{\bold q}}
	\def\MM{\mathcal M}
	\def\TT{\mathcal T}
	\def\EE{\mathcal E}
	\def\lsp{\mbox{lsp}}
	\def\rsp{\mbox{rsp}}
	\def\pf{\noindent {\it Proof.} }
	\def\mp{\mbox{pyramid}}
	\def\mb{\mbox{block}}
	\def\mc{\mbox{cross}}
	\def\qed{\hfill \rule{4pt}{7pt}}
	\def\block{\hfill \rule{5pt}{5pt}}
    \def\lr#1{\multicolumn{1}{|@{\hspace{.6ex}}c@{\hspace{.6ex}}|}{\raisebox{-.3ex}{$#1$}}}

	\begin{center}
		{\Large \bf  On Refinements of   Wilf-Equivalence for  Involutions}
	\end{center}
	
	\begin{center}
		{\small  Sherry H.F. Yan$^a$, Lintong Wang$^{a}$,  Robin D.P. Zhou$^{b,*}$\footnote{$^*$Corresponding author.}  \footnote{{\em E-mail address:} dapao2012@163.com.}}
		
		$^{a}$Department of Mathematics,
		Zhejiang Normal University\\
		Jinhua 321004, P.R. China
		
		$^b$College of Mathematics Physics and Information\\
		Shaoxing University\\
		Shaoxing 312000, P.R. China
		
	\end{center}
	
	\noindent {\bf Abstract.}
	Let $\mathcal{S}_n(\pi)$ (resp. $\mathcal{I}_n(\pi)$ and $\mathcal{AI}_n(\pi)$) denote the set of permutations (resp. involutions and  alternating involutions) of length $n$ which avoid the permutation pattern $\pi$. 	For $k,m\geq 1$,
Backelin-West-Xin  proved that $|\mathcal{S}_n(12\cdots k\tau)|=
|\mathcal{S}_n(k\cdots 21\tau)|$ by establishing a bijection between these two sets, where
$\tau = \tau_1\tau_2\cdots \tau_m$ is an arbitrary permutation of $k+1,k+2,\ldots,k+m$.
 The result has been extended to involutions by
Bousquet-M\'elou and Steingr\'imsson  and  to alternating permutations by the first author.
In this paper, we  shall  establish  a peak set preserving bijection between $\mathcal{I}_n(123\tau)$ and
$\mathcal{I}_n(321\tau)$  via transversals, matchings, oscillating tableaux and pairs of noncrossing Dyck paths as intermediate structures. Our result is a refinement of     the  result of  Bousquet-M\'elou  and Steingr\'imsson  for the case when $k=3$.
As an application,  we show bijectively  that $|\mathcal{AI}_n(123\tau)| =
|\mathcal{AI}_n(321\tau)|$, confirming a recent  conjecture of Barnabei-Bonetti-Castronuovo-Silimbani.
Furthmore, some conjectured equalities posed by Barnabei-Bonetti-Castronuovo-Silimbani  concerning pattern avoiding alternating
involutions  are also proved.
	


	%===========================================================================
	
		\section{Introduction}
	
	Let $[n]=\{1,2,\ldots, n\}$.
	Given a permutation $\pi=\pi_1 \pi_2 \cdots \pi_n$ of $[n]$, it is said to be {\em alternating} if $\pi_1 < \pi_2 > \pi_3 < \pi_4 > \cdots$
	and it is said
	to be {\em reverse alternating} if $\pi_1 > \pi_2 < \pi_3 > \pi_4 < \cdots$.
	A permutation $\pi$ is an {\em involution} if $\pi = \pi^{-1}$.
	Denote $\mathcal{S}_n$ (resp. $\mathcal{A}_n$,   $\mathcal{I}_n$, $\mathcal{AI}_n$ and $\mathcal{RAI}_n$) to be
	the set of permutations (resp. alternating permutations,   involutions, alternating involutions and reverse alternating involutions) of length $n$.
	
	\begin{example}
		The permutation $\pi = 4\,5\,3\,8\,1\,6\,2\,7$ is an alternating permutation but not an involution,
		while the permutation $\sigma = 6\, 4\, 8\, 2\, (10)\, 1\, 9\, 3\, 7\, 5$ is a reverse alternating involution.
	\end{example}

	In a permutation $\pi=\pi_1\pi_2
\cdots \pi_n$, an index    $i$ ($1<i<n$) is said to be a {\em peak} of $\pi$ if
	$\pi_{i-1} < \pi_i > \pi_{i+1}$.
	Let  $\mathrm{Peak}(\pi)$ denote the set of
	peaks of $\pi$. For example,  given $\pi = 5\,4\,7\,9\,8\,3\,6\,1\,2$,
		 we have	$\mathrm{Peak}(\pi) = \{4,7\}$.

	Given a permutation $\pi \in \mathcal{S}_n$ and a permutation $\sigma \in \mathcal{S}_k$,
	we say $\pi$   contains  the pattern $\sigma$ if there exists a subsequence of
	$\pi$ which is order isomorphic to $\sigma$.
	Otherwise, we say $\pi$ avoids the patter $\sigma$  and
	$\pi$ is $\sigma$-avoiding.
	Let $\mathcal{S}_n(\sigma)$ denote the set of  $\sigma$-avoiding permutations of length $n$.
	We will keep this notation also when $\mathcal{S}_n$ is replaced by
	other subsets of permutations such as $\mathcal{A}_n,\mathcal{I}_n,\mathcal{AI}_n,\mathcal{RAI}_n$, etc.
	
	\begin{example}
		For permutation $\pi = 5\,4\,7\,9\,8\,3\,6\,1\,2$,
		it is $1234$-avoiding while it contains the pattern
		$4321$.
	\end{example}

The notion of pattern avoiding permutations was introduced by Knuth \cite{Knuth} in 1970 and was first systematically   studied by Simion-Schmidt \cite{Simion}.
Pattern avoiding permutations have been extensively exploited  over particular subset of $\mathcal{S}_n$. For example, various results have been obtained for
	   pattern avoiding  alternating permutations
	(see e.g. \cite{Bona,Chen,Lewis2009,Lewis2011,Lewis2012,Stanley,Yan2012})  and pattern avoiding involutions (see e.g. \cite{Barnabei2011,Bona2016,Jaggard}).

Let $\alpha$ and $\beta$ be two permutations in $\mathcal{S}_k$ and
let $\mathcal{B}_n$ be a subset of $\mathcal{S}_n$.
We say that $\alpha$ and $\beta$ are {\em Wilf-equivalent} over the set $\mathcal{B}_n$ if $|\mathcal{B}_n(\alpha)| = |\mathcal{B}_n(\beta)|$ for all positive integers $n$. The {\em direct sum} of two permutations $\pi=\pi_1\pi_2\cdots \pi_k\in \mathcal{S}_k$  and $\tau=\tau_1\tau_2\cdots \tau_m\in \mathcal{S}_m$, denoted by $\pi\oplus \tau$, is the permutation $\pi_1\pi_2\cdots \pi_k (\tau_1+k)(\tau_2+k)\cdots (\tau_m+k)$.
We say that $\alpha$ and $\beta$ may be {\em prefix exchanged} if for any positive integer $n$, we have
  $|\mathcal{B}_n(\alpha\oplus\tau)|=|\mathcal{B}_n(\beta\oplus\tau)|$ holds for any nonempty permutation $\tau$.  In this context, the pair
 $(\alpha, \beta)$  is called  {\em  a prefix exchanging pair } over the set $\mathcal{B}_n$.
Note that the ability to exchange  two prefixes  would imply infinite class of Wilf-equivalent pairs.

Let $I_k=12\cdots k$ and $J_k=k\cdots 21$.  Most of the known results for prefix exchanging are related to the pair  $(I_k,J_k)$.    Backelin-West-Xin \cite{BWX} proved that
$|\mathcal{S}_n(I_k\oplus \tau)|=|\mathcal{S}_n(J_k\oplus \tau)|$ for general $k$ and any pattern $\tau$, which extends the results obtained by West \cite{West} for $k=2$ and  Babson-West $\cite{BW}$ for  $k=3$.

The systematic study of pattern avoiding involutions was also initiated in \cite{Simion},
and continued by Guibert \cite{Guibert} in his thesis for patterns of length $4$.
Guibert discovered that, for a large number of patterns $\tau$ of length $4$, $\mathcal{I}_n(\tau)$
is counted by the $n$-th Motzkin number:
$$M_n = \sum_{k=0}^{\lfloor n/2\rfloor}{n! \over k!(k+1)!(n-2k)!}.$$
In the spirit of the work of Babson-West \cite{BW},  Jaggard \cite{Jaggard}
 proved that
$I_k$ and $J_k$ are prefix exchanged over involutions $\mathcal{I}_n$ for $k = 2, 3$, and further conjectured that
  it also holds for general  $k$. This conjecture has been  confirmed by
Bousquet-M\'elou  and Steingr\'imsson \cite{Bousquet}.

Anaologious to ordinary permutations,
 pattern avoiding alternating permutations  have been extensively investigated, see    \cite{Bona,Chen,Lewis2009, Lewis2011,Lewis2012, Mansor2003, Stanley, Yan2012, Yan2013} and references therein.
  Yan \cite{Yan2013}
proved that $I_k$ and $J_k$ are  prefix exchanged for general $k$ over alternating permutations $\mathcal{A}_n$,  paralleling the  result of   Backelin-West-Xin \cite{BWX}  for ordinary permutations and the result of Bousquet-M\'elou  and Steingr\'imsson  \cite{Bousquet} for involutions.

	Recently, Barnabei-Bonetti-Castronuovo-Silimbani  \cite{Barnabei} initiated the study of pattern avoiding alternating involutions.
They proved that $I_2$ and $J_2$ are prefix exchanged over alternating permutations $\mathcal{AI}_n$.
	They also investigated the enumeration of  some classes of alternating and reverse alternating
	involutions avoiding a pattern of length $3$ or $4$. Interestingly, alternating   involutions and reverse alternating involutions avoiding a given pattern of length $4$ are showed to be closely related to Motzkin numbers in \cite{Barnabei}.
 Table \ref{table:Barnabei} lists the known enumeration  results  for alternating involutions avoiding a pattern of length $4$.
\begin{table}[!h]
\centering
\footnotesize
\caption{Alternating and reverse alternating involutions avoiding a pattern of length $4$}\label{table:Barnabei}
\vskip 2mm
\begin{tabular}{|c|c|} \hline
patterns & enumeration results \\ \hline
               &      $|\mathcal{RAI}_{2n}(4321)| = |\mathcal{AI}_{2n}(1234)|= M_n$,  \\[3pt]
$$1234$, 4321$ &    $|\mathcal{RAI}_{2n-1}(1234)| = |\mathcal{AI}_{2n-1}(4321)| = |\mathcal{RAI}_{2n-1}(4321)| =
                       |\mathcal{AI}_{2n-1}(1234)| =  M_n - M_{n-2}$,   \\[3pt]
               &    $|\mathcal{RAI}_{2n}(1234)| = |\mathcal{AI}_{2n}(4321)|= M_{n+1} - 2M_{n-1} + M_{n-3}$.  \\ \hline
$3412$         &   $|\mathcal{RAI}_{2n}(3412)| = |\mathcal{AI}_{2n+2}(3412)| = |\mathcal{AI}_{2n+1}(3412)| = |\mathcal{RAI}_{2n+1}(3412)| = M_n$.	\\ \hline
                &     $|\mathcal{AI}_{2n}(1243)| = |\mathcal{AI}_{2n}(2143)|= |\mathcal{AI}_{2n}(2134)| = M_n$, \\[3pt]
$2143$, $2134$  &  $|\mathcal{AI}_{2n-1}(2134)| = |\mathcal{RAI}_{2n-1}(1243)| =  M_n - M_{n-2}$, \\[3pt]
      $1243$    &  $|\mathcal{AI}_{2n+1}(1243)| = |\mathcal{AI}_{2n+1}(2143)|= |\mathcal{RAI}_{2n+1}(2143)| =   |\mathcal{RAI}_{2n+1}(2134)| = M_n$, \\                      &            $|\mathcal{RAI}_{2n}(2143)| = M_{n-1}$.  \\ \hline
\end{tabular}
\end{table}

In \cite{Barnabei}, Barnabei-Bonetti-Castronuovo-Silimbani  posed the following conjectures.
\begin{conjecture} \label{conj1}
	(\cite{Barnabei}, Conjecture 12.1)
	$$|\mathcal{RAI}_{2n}(1243)| = |\mathcal{RAI}_{2n}(2134)| = M_n.$$
\end{conjecture}

\begin{conjecture} \label{conj2}
	(\cite{Barnabei}, Conjecture 12.2)
	\begin{align}
		|\mathcal{AI}_{2n}(1432)| &= |\mathcal{AI}_{2n}(3214)|  = M_n,  \label{conj2_0} \\[3pt]
|\mathcal{AI}_{2n-1}(3214)| &= |\mathcal{RAI}_{2n-1}(1432)| = M_n - M_{n-2}, \label{conj2_3}\\[3pt]
|\mathcal{RAI}_{2n}(1432)| &= |\mathcal{RAI}_{2n}(3214)| = M_n, \label{conj2_1} \\[3pt]
		|\mathcal{AI}_{2n+1}(1432)| &= |\mathcal{RAI}_{2n+1}(3214)| = M_n. \label{conj2_2}
			\end{align}
\end{conjecture}

\begin{conjecture} \label{conj3}
	(\cite{Barnabei}, Conjecture 12.3)
	Let $n\geq 1$.
  $I_3$ and $J_3$ are prefix exchanged over alternating permutations $\mathcal{AI}_n$, that is,
the equality $
	|\mathcal{AI}_n(123 \oplus \tau)| =|\mathcal{AI}_n(321 \oplus \tau)|$  holds  for any   nonempty pattern $\tau$.
\end{conjecture}
As remarked by Barnabei-Bonetti-Castronuovo-Silimbani \cite{Barnabei},    the proofs of Conjectures  \ref{conj1} and \ref{conj2} would complete  the classification of the patterns of length $4$ with respect to Wilf-equivalence for alternating and reverse alternating involutions.

For a permutation
	$\pi=\pi_1\pi_2\cdots \pi_n \in \mathcal{S}_n$, its {\em reverse} $\pi^r$ is defined by
	$\pi^r_i=\pi_{n+1-i}$ for $1\leq i\leq n$.
	The {\em complement} of $\pi$, denoted by $\pi^c$, is defined by $\pi^c_i=n+1-\pi_i$ for  $1\leq i\leq n$.
	The {\em reverse-complement} of $\pi$, denoted by $\pi^{rc}$, is
	defined by $\pi^{rc}_i = n+1-\pi_{n+1-i}$.
	Notice that $(\pi^{rc})^{-1} = (\pi^{-1})^{rc}$.
	The reverse-complement preserves the property of being an involution.
	\begin{example}
		Let $\pi = 5\,4\,7\,9\,8\,3\,6\,1\,2$.
		We have $\pi^r = 2\,1\,6\,3\,8\,9\,7\,4\,5$,
		$\pi^c = 5\,6\,3\,1\,2\,7\,4\,9\,8$ and
		$\pi^{rc} = 8\,9\,4\,7\,2\,1\,3\,6\,5$.
	\end{example}	

	\begin{lemma}\label{lem:reverse_complement}
		The reverse-complement is a bijection between $\mathcal{AI}_{2n}(\tau)$ and
		$\mathcal{AI}_{2n}(\tau^{rc})$, between $\mathcal{RAI}_{2n}(\tau)$ and  $\mathcal{RAI}_{2n}(\tau^{rc})$ and between $\mathcal{AI}_{2n+1}(\tau)$ and
		$\mathcal{RAI}_{2n+1}(\tau^{rc})$.
	\end{lemma}

Combining Lemma  \ref{lem:reverse_complement} and Table \ref{table:Barnabei},  Conjecture \ref{conj3} would imply    (\ref{conj2_0}) and (\ref{conj2_3})  of Conjecture \ref{conj2} if we set $\tau = 1$ in Conjecture \ref{conj3}.
In order to prove Conjecture \ref{conj3},
   we shall derive    the following stronger result.

	\begin{theorem} \label{generalth}
		 		Let $n\geq 1$.  For any  pattern $\tau$, there exists a    bijection $\Phi$ between $\mathcal{I}_n(123 \oplus \tau)$ and $\mathcal{I}_{n}(321 \oplus \tau)$ such that for any $\pi\in \mathcal{I}_n(123 \oplus \tau)$, we have $\mathrm{Peak}(\pi)=\mathrm{Peak}(\Phi(\pi))$.
\end{theorem}

  Note that for any pattern $\tau$,   Bousquet-M\'elou  and Steingr\'imsson \cite{Bousquet}  showed that the equality   $|\mathcal{I}_n(I_k\oplus \tau)| =|\mathcal{I}_{n}(J_k \oplus \tau)|$ holds for   general $k$. Hence, Theorem \ref{generalth} can be
   viewed as a refinement of   the result of  Bousquet-M\'elou  and Steingr\'imsson \cite{Bousquet} for the case when $k=3$.
    We  remark that the bijection between  $\mathcal{I}_n(I_k\oplus \tau)$ and $\mathcal{I}_{n}(J_k \oplus \tau)$ established by  Bousquet-M\'elou  and Steingr\'imsson in \cite{Bousquet} can not
    preserve the set of peaks,
   while the bijection between  $\mathcal{A}_n(I_k\oplus \tau)$ and $\mathcal{A}_{n}(J_k \oplus \tau)$ established  by  Yan in \cite{Yan2013} can not
    preserve the property of being  involutions.

%Actually, we will be proving an extended version of the above theorem as shown in Theorem \ref{mainth1}.
%We extend the definition of pattern avoidance to
%where we extend the
%Actually, we extend the result to cover Young diagram.
%And the above theorem follows from the  prefix exchanging for shape-Wilf-equivalence given in
%Theorem \ref{mainth1}.
%We should note that
%the above theorem is actually  deduce from an extended version as shown in Theorem \ref{mainth1},
%where we extend the definition of pattern avoidance to a transversal of a Young diagram  defined in
%Section \ref{sec:3}.

%The objective of this paper is to prove the above conjectures.
This paper is organized as follows.
In Section \ref{sec:2}, we shall prove Conjecture \ref{conj1} by building
a bijection between $\mathcal{RAI}_{2n}(1243)$ and a subset of $\mathcal{RAI}_{2n}(1234)$, the latter of which will be showed to be counted by $M_n$.
In Section \ref{sec:3}, we shall establish the bijection $\Phi$ as stated in Theorem \ref{generalth} via transversals, matchings, oscillating tableaux and pairs of noncrossing Dyck paths as intermediate structures. Furthermore, we  prove  that  $\Phi$ induces a bijection between $\mathcal{AI}_n(321 \oplus \tau)$ and $\mathcal{AI}_n(123 \oplus \tau)$ for any nonempty pattern $\tau$, thereby confirming Conjecture \ref{conj3}.

	\section{Proof of Conjecture \ref{conj1}} \label{sec:2}
	
	In this section, we construct a bijection between $\mathcal{RAI}_{2n}(1243)$ and
	$\{\pi \in \mathcal{RAI}_{2n}(1234) \mid \pi_1 = 2n \, \text{or} \, \pi_3 = 2n\}.$
	We prove that  the cardinality of the latter set is given by $M_n$ by a series of refined decompositions,
    thereby confirming Conjecture \ref{conj1}.

   Given a permutation $\pi \in \mathcal{S}_n$, a {\em descent} in a permutation $\pi$ is an index $i$ ($1 \leq i < n$) such that 	$\pi_i > \pi_{i+1}$.
	An {\em ascent} in a permutation $\pi$ is an index $i$ ($1 \leq i < n$) such that
	$\pi_i < \pi_{i+1}$.
	Let $\mathrm{Des}(\pi)$ and $\mathrm{Asc}(\pi)$ denote the set of
	descents and ascents of $\pi$, respectively.
\begin{example}
  Let $\pi = 5\,4\,7\,9\,8\,3\,6\,1\,2$.
		We have $\mathrm{Des}(\pi) = \{1,4,5,7\}$,
		$\mathrm{Asc}(\pi) = \{2,3,6,8\}$.
\end{example}

	Let us recall a bijection $f$ between $\mathcal{S}_n(I_k)$ and
	$\mathcal{S}_n(I_{k-2}\oplus 21)$ which was first defined by West \cite{West}.
	Given a permutation $\pi = \pi_1\pi_2\cdots\pi_n \in \mathcal{S}_n$,
	the {\em rank} of the element $\pi_i$ is defined to be the maximum length of
	an increasing subsequence ending at $\pi_i$.
	For $\pi\in  \mathcal{S}_n(I_k)$, the maximal rank of $\pi$
	is $k-1$.
	Let $\mathcal{E}$ denote the set of elements in $\pi$ of rank $k-1$ and
	let $\mathcal{L}$ denote the set of their positions.
	Then $\sigma = \sigma_1\sigma_2\cdots \sigma_n = f(\pi)$ is defined as follows.
	\begin{itemize}
		\item If $i \notin \mathcal{L}$, $\sigma_i = \pi_i$.
		\item If $i \in \mathcal{L}$, $\sigma_i$ is the smallest unused element of $\mathcal{E}$
		which is larger that the closest entry of rank $k-2$ to the left of $\pi_i$.
	\end{itemize}
	
	\begin{example}
		Given $\pi = 6\, 4\, 8\, 2\, (10)\, 1\, 9\, 3\, 7\, 5 \in \mathcal{S}_{10}(1234)$, the elements of rank $1$ are $6,4,2$ and $1$,  the elements of rank $2$ are $8$
		and $3$, and the elements of rank $3$ are $10,9,7$ and $5$. Hence, we have $\mathcal{E}=\{10,9,7,5\}$ and $\mathcal{L}=\{5, 7, 9, 10\}$.
		Then we have $f(\pi) = 6\, 4\, 8\, 2\, 9\, 1\, (10)\, 3\, 5\, 7$.
	\end{example}

	From the definition of $f$, we can see that $f$ fixes the elements of
	rank $k-2$ or less.
	Given a permutation $\sigma = \sigma_1\sigma_2\cdots\sigma_n \in \mathcal{S}_n(I_k\oplus 21)$,
	let $\mathcal{E}$ denote the set of elements in $\sigma$ with rank $k-1$ or more and
	let $\mathcal{L}$ denote the set of their positions.
	Then $\pi = \pi_1\pi_2\cdots \pi_n = f^{-1}(\sigma)$ is defined as follows.
	\begin{itemize}
		\item If $i \notin \mathcal{L}$, $\pi_i = \sigma_i$.
		\item If $i \in \mathcal{L}$, $\pi_i$ is the largest unused element of $\mathcal{E}$.
	\end{itemize}
	In other words, $\pi$ is obtained from $\sigma$ by rearranging the
	elements of rank $k-1$ or more of $\sigma$ in decreasing order and keeping elements of rank $k-2$ or less unchanged.
	It can be easily checked that $f^{-1}$ maps the elements of rank $k-1$ or
	more in $\sigma$ to the elements of rank $k-1$ in $\pi$.
	
	\begin{example}
		Given $\sigma = 6\, 4\, 8\, 2\, 9\, 1\, (10)\, 3\, 5\, 7 \in \mathcal{S}_{10}(1243)$, the elements of rank $1$ are $6,4,2$ and $1$, the elements of rank $2$ are $8$ and $3$, and the elements of rank $3$ or more are $9,10,5$ and $7$.
Hence, we have $\mathcal{E}=\{9,10,5,7\}$ and $\mathcal{L}=\{5, 7, 9, 10\}$.
		Then $f^{-1}(\sigma) = 6\, 4\, 8\, 2\, (10)\, 1\, 9\, 3\, 7\, 5$.
	\end{example}

	B\'ona \cite{Bona} proved that $f$ preserves the alternating
	property when the length of the permutation is even and used it
	to prove a generalized version of a conjecture of Lewis \cite{Lewis2012}.
	Barnabei-Bonetti-Castronuovo-Silimbani \cite{Barnabei} showed that $f$ preserves the
	property of being an involution and obtained some enumeration result
	about pattern avoiding alternating involutions.

    Given a word $w = w_1w_2\cdots w_k$ with distinct numbers, the {\em standardization} of $w$ is the unique permutation in $S_k$ which is
	order isomorphic to $w$.

It is routine to check that Conjecture \ref{conj1} holds for $n \leq 3$.
In the rest of this section, unless otherwise specified, we always assume that $n \geq 4$.
	In the following lemma, we study the case when $f$ restricts to the reverse alternating involutions in  $\mathcal{RAI}_{2n}(1243)$.
	
	\begin{lemma} \label{lem:f}
		 The map $f^{-1}$ is a bijection between $\mathcal{RAI}_{2n}(1243)$ and
		$\{\pi \in \mathcal{RAI}_{2n}(1234) \mid \pi_1 = 2n \; \text{or} \; \pi_3 = 2n\}$.
	\end{lemma}
\pf
Let $\sigma = \sigma_1\sigma_2 \cdots \sigma_{2n} \in \mathcal{RAI}_{2n}(1243)$ and let $\pi =\pi_1\pi_2 \cdots \pi_{2n}= f^{-1}(\sigma)$.
First we assert that $\sigma_1 = 2n$ or $\sigma_3 = 2n$.
Assume to contrary that $2n$ is not in position $1$ or $3$.
Since $\sigma$ is an involution and $\sigma_{2n-1} > \sigma_{2n}$, we have that $2n$ appears to the left of $2n-1$ in $\sigma$.
Then $\sigma_2\sigma_3(2n)(2n-1)$ would  form  a $1243$ pattern, a contradiction.
Hence the assertion is proved.
Then, it is easily seen that $2n$ is of rank $1$ or $2$ in $\sigma$.
%As $\sigma$ is an involution,  $\sigma_{2n}$ is also of rank $1$ or $2$.
Since $f^{-1}$ fixes the elements of rank $1$ and $2$, we have
$\pi_1 = 2n$ or $\pi_3 = 2n$.

Since $f^{-1}$ is a bijection between $\mathcal{I}_{2n}(1243)$ and
$\mathcal{I}_{2n}(1234)$,
to prove that $\pi \in \mathcal{RAI}_{2n}(1234)$,
it suffices to show that $\pi$ is reverse alternating.  More precisely, we aim to show that $\mathrm{Peak}(\pi)=\mathrm{Peak}(\sigma)$ and $\pi_1>\pi_2$.
Notice that $f^{-1}$ fixes the elements of rank $1$ or $2$ in $\sigma$
and maps the elements of rank $3$ or
	more in $\sigma$ to the elements of rank $3$ in $\pi$. Hence, we have $\pi_1=\sigma_1>\sigma_2=\pi_2$ as both  $\sigma_1$ and $\sigma_2$ are of rank $1$.
Let $\sigma_k$ be an element in $\sigma$ of rank $3$ or more.
Then there exist $i<j<k$ such that  $\sigma_i\sigma_j\sigma_k$  forms a
$123$ pattern of $\sigma$.
We assert that  $k\in \mathrm{Peak}(\sigma)$.
If not, then $\sigma_i\sigma_j\sigma_{k-1}\sigma_k$ would
form a $1243$ pattern, a contradiction. By similar arguments, one can verify that  both $\sigma_{k-1}$  and $\sigma_{k+1}$ are of rank $2$ or less.
As $f^{-1}$ fixes the elements of rank $2$ or less in $\sigma$  and maps the elements of rank $3$ or
more in $\sigma$ to the elements of rank $3$ in $\pi$,
then both $\pi_{k-1}$ and $\pi_{k+1}$ are of rank $2$ or less,  and $\pi_k$ is an element of rank $3$ in $\pi$.
This implies that $\pi_{k-1} < \pi_k >\pi_{k+1}$, that is,
 $k \in \mathrm{Peak}(\pi)$. So far, we have  concluded that the positions of the elements of rank $3$ or more in $\sigma$ are peaks of $\sigma$, and $f^{-1}$ preserves the property of being a peak. Therefore,   we have $\mathrm{Peak}(\pi)=\mathrm{Peak}(\sigma)$ as desired.

To prove that $f^{-1}$ is a bijection between $\mathcal{RAI}_{2n}(1243)$ and
$\{\pi \in \mathcal{RAI}_{2n}(1234) \mid \pi_1 = 2n \; \text{or} \; \pi_3 = 2n\}$,
it remains to show that given $\pi \in \{\pi \in \mathcal{RAI}_{2n}(1234) \mid \pi_1 = 2n \; \text{or} \; \pi_3 = 2n\}$, $f(\pi)$ is a permutation in $\mathcal{RAI}_{2n}(1243)$.
Since $f$ is a bijection between $\mathcal{I}_{2n}(1234)$ and $\mathcal{I}_{2n}(1243)$,
it is sufficient to show that $f(\pi)$ is reverse alternating.
By similar arguments as above, the positions of the elements of rank $3$ in $\pi$ are peaks of $\pi$, and $f$ preserves the property of being a peak.
It follows that $f(\pi) \in \mathcal{RAI}_{2n}(1243)$,
completing the proof.
\qed

	In the following, we proceed to prove  $|\{\pi \in \mathcal{RAI}_{2n}(1234) \mid \pi_1 = 2n \; \text{or} \; \pi_3 = 2n\}| = M_n$  by distinguishing the value of $\pi_1$.
	Given $\pi \in \mathcal{RAI}_{2n}(1234)$ with $\pi_1 = 2n$,
	the permutation $\pi'$ obtained from $\pi$ by removing
	the elements $2n$ and $1$ is an arbitrary permutation in $\mathcal{AI}_{2n-2}(1234)$, which is counted by $M_{n-1}$ according to Table \ref{table:Barnabei}.
Therefore, we deduce that
\begin{equation}\label{equ:p1=1}
  |\{\pi \in \mathcal{RAI}_{2n}(1234) \mid \pi_1 = 2n\}| = M_{n-1}.
\end{equation}

	By Lemma \ref{lem:reverse_complement} and Lemma \ref{lem:f}, to prove Conjecture \ref{conj1},
	it suffices to prove $|\{\pi \in \mathcal{RAI}_{2n}(1234) \mid \pi_3 = 2n \}| = M_n - M_{n-1}$.
Let $\mathcal{O}_{2n} = \{\pi \in \mathcal{RAI}_{2n}(1234) \mid \pi_3 = 2n \}$.
	Given a permutation $\pi \in \mathcal{O}_{2n}$, if we delete the elements $2n$ and $3$ from $\pi$ and standardize, we obtain an
	arbitrary permutation in $\mathcal{P}_{2n-2}\bigcup \mathcal{Q}_{2n-2}$, where

		$$\mathcal{P}_{2n-2} = \{\pi \in \mathcal{I}_{2n}(1234)\mid \mathrm{Des}(\pi) = \{1\}\cup\{2,4,6,\ldots, 2n-4 \}\}$$
and
		$$\mathcal{Q}_{2n-2}  = \{\pi \in \mathcal{I}_{2n}(1234)\mid \mathrm{Des}(\pi) = \{1\}\cup\{4,6,8,\ldots,2n-4 \}\}.$$
	
%	We proceed to show that $|\mathcal{P}_{2n}| = M_{n+1} - M_n - M_{n-1}$ and
%	$|\mathcal{Q}_{2n}| = M_{n-1}$.
    Next we will be devoted to the enumeration of $\mathcal{P}_{2n}$ and $\mathcal{Q}_{2n}$.
	The following lemma is needed for the enumeration of $\mathcal{P}_{2n}$.
	
	\begin{lemma}\label{lem:R}
	Let $\mathcal{R}_{2n+1} = \{\pi \in \mathcal{I}_{2n+1}(1234) \mid \mathrm{Des}(\pi) = \{3,5,7,\ldots, 2n-1\}\cup\{2n\}\}$ and
		$\mathcal{R}_{2n+1}^{rc} = \{\pi \in \mathcal{I}_{2n+1}(1234) \mid \mathrm{Des}(\pi) = \{1\}\cup\{2,4,6,\ldots, 2n-2\}\}$.
		Then we have
		\begin{align*}
			|\mathcal{R}_{2n+1}| = |\mathcal{R}_{2n+1}^{rc}| =  M_{n} - M_{n-1} - M_{n-2}.
		\end{align*}
	\end{lemma}
	
	\pf
	Notice that the map reverse-complement is a bijection between $\mathcal{R}_{2n+1}$ and $\mathcal{R}_{2n+1}^{rc}$.
  To prove this lemma, it suffices to show
	 $|\mathcal{R}_{2n+1}| = M_{n} - M_{n-1} - M_{n-2}$.
  %It can be easily checked that the conclusion is true for $n \leq 1$.
	%In the following, unless otherwise specified,
   Let $\pi$ be a permutation in $\mathcal{R}_{2n+1}$.
    From the fact that $\pi$ is a $1234$-avoiding involution, we must have $\pi_3 = 2n+1$ and
    $\pi_{2n+1} = 3$. If not, $\pi_1\pi_2\pi_3(2n+1)$ would form a $1234$ pattern, a contradiction.
	Moreover, we claim that $\pi_2>\pi_4$.
    Otherwise, $\pi_1\pi_2\pi_4\pi_5$ would form a $1234$ pattern of $\pi$, a contradiction.
	Let $\tau$ be the permutation obtained from $\pi$ by removing the elements
	$2n+1$ and $3$ and standardizing.
   Since $\pi_{2n} > \pi_{2n+1} = 3$, it follows that
   $\tau_{2n-1} = \pi_{2n} -1 > 2$.
	Then $\tau$ is an arbitrary permutation in
	$\{\pi \in \mathcal{AI}_{2n-1}(1234) \mid \pi_2 \neq 2n-1\}$.
	By Table \ref{table:Barnabei}, we have $|\mathcal{AI}_{2n-1}(1234)| = M_{n} - M_{n-2}$.
	Since $|\{\pi \in \mathcal{AI}_{2n-1}(1234) \mid \pi_2 = 2n-1\}| = M_{n-1}$ (See \cite{Barnabei} Theorem 7.1), we derive that
	\begin{align*}
		|\{\pi \in \mathcal{AI}_{2n-1}(1234) \mid \pi_2 \neq 2n-1\}| &= |\mathcal{AI}_{2n-1}(1234)| - |\{\pi \in \mathcal{AI}_{2n-1}(1234) \mid \pi_2 = 2n-1\}|  \\[3pt]
		&= M_n - M_{n-2} - M_{n-1},
	\end{align*}
	as desired, completing the proof.
	\qed
	
	Let us recall some notation and terminology.
	Given an integer $n$, a {\em partition} of $n$ is a sequence of
	nonnegative integers $\lambda=(\lambda_1,\lambda_2,\ldots,\lambda_k)$
	such that $n=\lambda_1+\lambda_2+\cdots+\lambda_k$ and $\lambda_1 \geq \lambda_2 \geq \cdots \geq \lambda_k >0$.
	We denote $\lambda \vdash n$.
	Each $\lambda_i$ is called a {\em part} of $\lambda$.
	A {\em Young diagram} of shape $\lambda$ is defined to be a left-justified array of $n$ boxes with $\lambda_1$ boxes
	in the first row, $\lambda_2$ boxes in the second row and so on.
	%We call $n$ the size of the Young diagram.
	
	A standard Young tableau (SYT) $T$ of shape $\lambda$ is a filling
	of the Young diagram of shape $\lambda$ with the numbers
	$1,2,\ldots,n$ such that the entries are increasing alone each rows and
	each columns.
	Figure \ref{fig:SYT} illustrates an SYT of
	shape $(4,3,2,2)$.
	
	\begin{figure} [h]
		\centering
			\raisebox{-.6ex}{$\begin{array}[b]{*{4}c}\cline{1-4}
					\lr{1}&\lr{2}&\lr{5}&\lr{8}\\\cline{1-4}
					\lr{3}&\lr{4}&\lr{7}\\\cline{1-3}
					\lr{6}&\lr{9}\\\cline{1-2}
					\lr{10}&\lr{11}\\\cline{1-2}
				\end{array}$}
		\caption{An example of an SYT of shape $(4,3,2,2)$.}\label{fig:SYT}
	\end{figure}
	
	Given an SYT $T$, a {\em descent} of $T$
	is an element $k$ such that $k+1$ is below $k$ (not necessary directly below) and an {\em ascent} of $T$
	is an element $k$ such that $k+1$ is on the right of $k$ (not necessary directly right).
	Denote by $\mathrm{Des}(T)$ and $\mathrm{Asc}(T)$ the set of
	descents and ascents of $T$, respectively.
    For example, let $T$ be an SYT as shown in Figure \ref{fig:SYT}.
    Then we have $\mathrm{Des}(T) = \{2,5,8,9\}$ and $\mathrm{Asc}(T) = \{1,3,4,6,7,10\}$.
	
	The Robinson-Schensted-Knuth map (RSK) \cite{StanleyVol2} builds a bijection between
	symmetric group $\mathcal{S}_n$ and pairs $(P,Q)$ of standard Young tableaux of
	the same shape $\lambda \vdash n$.
	We denote this correspondence by $\pi \xrightarrow{\mathrm{RSK}} (P,Q)$, where $\pi \in \mathcal{S}_n$.
	By the property of RSK algorithm,
	the length of the longest increasing (resp. decreasing) subsequence of $\pi$ is equal to the number of columns (resp. rows) of $P$.
	Moreover, it can be easily seen that $\mathrm{Des}(\pi) = \mathrm{Des}(Q)$ by the insertion rule of
RSK algorithm.
Notice that  $\pi \xrightarrow{\mathrm{RSK}} (P,Q)$ if only if $\pi^{-1}\xrightarrow{\mathrm{RSK}} (Q,P)$.
	It follows that RSK algorithm associates an involution $\pi$
	with an SYT $T$ such that $\mathrm{Des}(\pi) = \mathrm{Des}(T)$.
    We also denote this correspondence by  $\pi \xrightarrow{\mathrm{RSK}} T$ if $\pi$ is an involution.
    Let $\mathrm{RSK}^{-1}$ denote the inverse of the RSK algorithm.
    We will simply write
    $T \xrightarrow{\mathrm{RSK}^{-1}} \pi$ if $\pi$ corresponds to the pair $(T,T)$ under the map $\mathrm{RSK}^{-1}$.
	
	\begin{lemma}\label{lem:P}
		There exists a bijection between $\mathcal{P}_{2n-2}$ and
		$\mathcal{R}_{2n+1}^{rc}$.
	\end{lemma}
	
	\pf
	Recall that
	$$\mathcal{P}_{2n-2} = \{\pi \in \mathcal{I}_{2n-2}(1234)\mid \mathrm{Des}(\pi) = \{1\}\cup\{2,4,6,\ldots, 2n-4 \}\}$$ and
	$$\mathcal{R}_{2n+1}^{rc} = \{\pi \in \mathcal{I}_{2n+1}(1234) \mid \mathrm{Des}(\pi) = \{1\}\cup\{2,4,6,\ldots, 2n-2\} \}.$$
	By the above analysis, RSK algorithm builds a bijection
	between $\mathcal{I}_n(1234)$ and standard Young tableaux with at most
	three columns which preserves the  set of	descents.

	Given $\pi \in \mathcal{P}_{2n-2}$, we construct a map $\gamma: \mathcal{P}_{2n-2}\rightarrow \mathcal{R}_{2n+1}^{rc}$
	as follows:
	\begin{equation*}
		\pi \xrightarrow{\mathrm{RSK}} T \xrightarrow{\eta} T' \xrightarrow{\mathrm{RSK}^{-1}} \sigma,
	\end{equation*}
	where $T' = \eta(T)$ is the SYT obtained
	from $T$ by adding $2n-1$ (resp. $2n$ and $2n+1$) in the bottom of the first (resp. second and third) column of $T$. See Figure \ref{fig:gamma} as an example.

\begin{figure}[H]
  \centering
  $\pi = 3\, 2\, 1\, 5\, 4\, 6 \xrightarrow{\mathrm{RSK}} T=\begin{array}[c]{*{3}c}\cline{1-3}
\lr{1}&\lr{4}&\lr{6}\\\cline{1-3}
\lr{2}&\lr{5}\\\cline{1-2}
\lr{3}\\\cline{1-1}
\end{array}
\xrightarrow{\enspace\eta\enspace}
    T'=\begin{array}[c]{*{3}c}\cline{1-3}
\lr{1}&\lr{4}&\lr{6}\\\cline{1-3}
\lr{2}&\lr{5}&\lr{9}\\\cline{1-3}
\lr{3}&\lr{8}\\\cline{1-2}
\lr{7}\\\cline{1-1}
\end{array}
\xrightarrow{\mathrm{RSK}^{-1}} \sigma = 7\, 3\, 2\, 8\, 5\, 9\, 1\, 4\, 6$.
  \caption{An example for the bijection $\gamma$.}\label{fig:gamma}
\end{figure}

	First we need to show that $\sigma = \gamma(\pi)$ is a permutation in $\mathcal{R}_{2n+1}^{rc}$.
Since $\pi$ is $1234$-avoiding, $T$ has at most three columns.
	By the definition of $\eta$, $T'$ contains three columns.
	If follows that $\sigma$ is $1234$-avoiding.
	As $2n-1$ is in the bottom of the first column of $T'$, we have
	$2n-1$ is below $2n-2$.
	Hence $2n-2 \in \mathrm{Des}(T')$.
	Since $2n$ is on the right of $2n-1$ and $2n+1$ is
	on the right of $2n$,  $2n-1$ and $2n$ are ascents of $T'$.
	Hence we have $\mathrm{Des}(T') = \mathrm{Des}(T) \cup \{2n-2\}$.
	Then
	$$\mathrm{Des}(\sigma) = \mathrm{Des}(T') = \mathrm{Des}(T) \cup \{2n-2\} =\mathrm{Des}(\pi) \cup \{2n-2\}= 	\{1\} \cup \{2,4,6,\ldots, 2n-2\},$$
	as RSK algorithm preserves the set of descents.
	To summarize, we have $\sigma = \gamma(\pi) \in \mathcal{R}_{2n+1}^{rc}$, the map
	$\gamma$ is well defined.
	
	We proceed to show that $\gamma$ is a bijection.
From the definition of $\gamma$, it suffices to show that
$\eta$ is a invertible.
	Assume that $\sigma \in R_{2n+1}^{rc}$ and let  $\sigma \xrightarrow{\mathrm{RSK}} T'$.
	Since $\sigma_{2n-1} < \sigma_{2n} < \sigma_{2n+1}$,
	we have $2n$ is on the right of $2n-1$, and
	$2n+1$ is on the right of $2n$ in $T'$.
   Notice that $\sigma$ is $1234$-avoiding, we derive that
   $T'$ contains at most three columns.
	It follows that $2n-1$ (resp. $2n$ and $2n+1$) is in the bottom of the first
	(resp. second and third) column of $T'$.
	Thus $\eta$ is invertible, completing the proof.
	\qed

Combining Lemma \ref{lem:R} and Lemma \ref{lem:P}, we obtain that
\begin{equation}\label{equ:P}
  |\mathcal{P}_{2n-2}| = M_n - M_{n-1} - M_{n-2}.
\end{equation}
	Now we proceed to enumerate the set $\mathcal{Q}_{2n}$.
	Recall that
	$$\mathcal{Q}_{2n} = \{\pi \in \mathcal{I}_{2n}(1234)\mid \mathrm{Des}(\pi) = \{1\}\cup\{4,6,8,\ldots,2n-2 \}\}.$$
	
	\begin{lemma} \label{lem:Q}
		We have that  $|\mathcal{Q}_{2n}| = M_{n-1}$ for $n\geq 3$.
	\end{lemma}
	
	\pf
	We would prove this lemma by induction on $n$.
	It is routine to check that the assertion is true
	for $n =3,4$.
	Assume that $n \geq 5$ and  the assertion holds for $  k < n$.
	
	Given $\pi \in \mathcal{Q}_{2n}$, we have that $2n$ appears to the right of $2n-1$ in $\pi$ since $\pi_{2n} > \pi_{2n-1}$.
	It can be deduced that $\pi_4 = 2n$.
	Otherwise, $\pi_2\pi_3\pi_4 (2n)$ would form a $1234$ pattern, a contradiction.
	From the fact that $\pi$ is an involution, we have
	$\pi_{2n} = 4$.
    This implies that $\pi_{2n-1} \leq 3$ as $\pi_{2n} > \pi_{2n-1}$.
    Actually, $\pi_{2n-1} = 1 \;\mbox{or} \; 3$.
    If not, we have $\pi_2 = 2n-1$ implying that
     $\pi_1 < \pi_2$. This yields a contradiction with the fact that  $1 \in \mathrm{Des}(\pi)$.
	Since $\pi$ is $1234$-avoiding, we have $\pi_3 > \pi_5$.
	If not, $\pi_2\pi_3\pi_5\pi_6$ would form a $1234$ pattern of $\pi$, a contradiction.
	We proceed to enumerate $\mathcal{Q}_{2n}$ by distinguishing the value of
	$\pi_{2n-1}$. \\
	\noindent {\bf Case  (\rmnum{1})}  $\pi_{2n-1} = 1$.
	By deleting the elements $2n-1,2n,1$ and $4$ from $\pi$ and
	standardizing,
	we obtain an arbitrary permutation in $\mathcal{AI}_{2n-4}(1234)$,
	which is counted by $M_{n-2}$ according to Table \ref{table:Barnabei}.

	\noindent {\bf Case  (\rmnum{2})}  $\pi_{2n-1} = 3$.
	By deleting the elements $2n-1$ and $3$ from $\pi$ and standardizing,
	we obtain an arbitrary permutation in $\mathcal{O}_{2n-2}$.
%It is routine to check that $|\mathcal{O}_{6}|=2=M_3-M_2$.  For $n\geq 5$,
	By (\ref{equ:P}) and  induction hypothesis,   we have
	\begin{equation} \label{equ:O}
		\begin{aligned}
			|\mathcal{O}_{2n-2}| &= |\{\pi \in \mathcal{RAI}_{2n-2}(1234) \mid \pi_3 = 2n-2 \}| = |\mathcal{P}_{2n-4}| + |\mathcal{Q}_{2n-4}|  \\
			&= M_{n-1} - M_{n-2} - M_{n-3} + M_{n-3} = M_{n-1} - M_{n-2}.
		\end{aligned}
	\end{equation}

	Combining Cases (\rmnum{1}) and (\rmnum{2}), we have reached the conclusion that
	$|\mathcal{Q}_{2n}| = M_{n-1}$  as desired.
	\qed
	
	Now we are ready for the proof of Conjecture \ref{conj1}.

	\noindent{\bf Proof of Conjecture \ref{conj1}.}
%    It is routine to check that the assertion holds for $n \leq 3$.
%    Now we assume that $n \geq 4$.
	By Lemma \ref{lem:reverse_complement} and  Lemma \ref{lem:f}, it is sufficient to prove that
	$$|\{\pi \in \mathcal{RAI}_{2n}(1234) \mid \pi_1 = 2n \;\text{or}\; \pi_3 = 2n\}| = M_n.$$
	By (\ref{equ:p1=1}), we have $|\{\pi \in \mathcal{RAI}_{2n}(1234) \mid \pi_1 = 2n \}| = M_{n-1}$.
By (\ref{equ:P}) and Lemma \ref{lem:Q},
	we have
	\begin{align*}
		 |\{\pi \in \mathcal{RAI}_{2n}(1234) \mid \pi_3 = 2n \}| &= |\mathcal{P}_{2n-2}| + |\mathcal{Q}_{2n-2}|   = M_n - M_{n-1}.
	\end{align*}
This completes the proof of Conjecture \ref{conj1}.
	\qed

	\section{Proof of  Theorem  \ref{generalth}}\label{sec:3}
	
	Let us begin with some necessary definitions and notations.  In a Young
	diagram, the square $(i,j)$ is referred to the square in column $i$ and row $j$, where columns are numbered from left to right and  rows are numbered  from top to bottom.   The {\em conjugate} of a partition $\lambda$, denoted by $\lambda^{T}$, is the partition whose Young diagram is the reflection along the main diagonal of $\lambda$'s Young diagram, and $\lambda$ is said to be {\em self-conjugate} if   $\lambda=\lambda^{T}$.
	
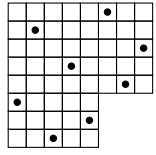
\begin{figure}[H]
\begin{center}
	\begin{tikzpicture}[font =\small , scale = 0.4]
		\draw (5,0)rectangle(5.6,0.6);
		\draw (5,0.6)rectangle(5.6,1.2);
		\draw (5,1.2)rectangle(5.6,1.8);
		\draw (5,1.8)rectangle(5.6,2.4);
		\draw (5,2.4)rectangle(5.6,3);
		\draw (5,3)rectangle(5.6,3.6);
		\draw (5,3.6)rectangle(5.6,4.2);
		\draw (5,4.2)rectangle(5.6,4.8);
		\draw (5.6,0)rectangle(6.2,0.6);
		\draw (5.6,0.6)rectangle(6.2,1.2);
		\draw (5.6,1.2)rectangle(6.2,1.8);
		\draw (5.6,1.8)rectangle(6.2,2.4);
		\draw (5.6,2.4)rectangle(6.2,3);
		\draw (5.6,3)rectangle(6.2,3.6);
		\draw (5.6,3.6)rectangle(6.2,4.2);
		\draw (5.6,4.2)rectangle(6.2,4.8);
		\draw (6.2,0)rectangle(6.8,0.6);
		\draw (6.2,0.6)rectangle(6.8,1.2);
		\draw (6.2,1.2)rectangle(6.8,1.8);
		\draw (6.2,1.8)rectangle(6.8,2.4);
		\draw (6.2,2.4)rectangle(6.8,3);
		\draw (6.2,3)rectangle(6.8,3.6);
		\draw (6.2,3.6)rectangle(6.8,4.2);
		\draw (6.2,4.2)rectangle(6.8,4.8);
    	\draw (6.8,0)rectangle(7.4,0.6);
		\draw (6.8,0.6)rectangle(7.4,1.2);
		\draw (6.8,1.2)rectangle(7.4,1.8);
		\draw (6.8,1.8)rectangle(7.4,2.4);
		\draw (6.8,2.4)rectangle(7.4,3);
		\draw (6.8,3)rectangle(7.4,3.6);
		\draw (6.8,3.6)rectangle(7.4,4.2);
		\draw (6.8,4.2)rectangle(7.4,4.8);
		\draw (7.4,0)rectangle(8,0.6);
		\draw (7.4,0.6)rectangle(8,1.2);
		\draw (7.4,1.2)rectangle(8,1.8);
 	    \draw (7.4,1.8)rectangle(8,2.4);
		\draw (7.4,2.4)rectangle(8,3);
		\draw (7.4,3)rectangle(8,3.6);
		\draw (7.4,3.6)rectangle(8,4.2);
		\draw (7.4,4.2)rectangle(8,4.8);
		\draw (8,1.8)rectangle(8.6,2.4);
		\draw (8,2.4)rectangle(8.6,3);
		\draw (8,3)rectangle(8.6,3.6);
		\draw (8,3.6)rectangle(8.6,4.2);
		\draw (8,4.2)rectangle(8.6,4.8);
		\draw (8.6,1.8)rectangle(9.2,2.4);
		\draw (8.6,2.4)rectangle(9.2,3);
		\draw (8.6,3)rectangle(9.2,3.6);
		\draw (8.6,3.6)rectangle(9.2,4.2);
		\draw (8.6,4.2)rectangle(9.2,4.8);
		\draw (9.2,1.8)rectangle(9.8,2.4);
		\draw (9.2,2.4)rectangle(9.8,3);
		\draw (9.2,3)rectangle(9.8,3.6);
		\draw (9.2,3.6)rectangle(9.8,4.2);
		\draw (9.2,4.2)rectangle(9.8,4.8);
		\filldraw[black](5.3,1.5)circle(3pt);
		\filldraw[black](5.9,3.9)circle(3pt);
		\filldraw[black](6.5,0.3)circle(3pt);
		\filldraw[black](7.1,2.7)circle(3pt);
		\filldraw[black](7.7,0.9)circle(3pt);
		\filldraw[black](8.3,4.5)circle(3pt);
		\filldraw[black](8.9,2.1)circle(3pt);
		\filldraw[black](9.5,3.3)circle(3pt);
				
	\end{tikzpicture}
\end{center}
\caption{A transversal of the Young diagram $\lambda=(8,8,8,8,8,5,5,5)$.}\label{transversal}
\end{figure}

	A  {\em $01$-filling} of a Young diagram $\lambda$   is obtained by filling   the squares of $\lambda$ with $1's$ and $0's$.
	A {\em  transversal} of  a Young diagram $\lambda$  is a  $01$-filling  of $\lambda$ such that every row and column contains  exactly one $1$, see Figure \ref{transversal} for  an  example,  where we represent a $1$ by a $\bullet$ and suppress the $0$'s.
Let $\mathcal{T}_n$ denote the set of transversals of all Young diagrams $\lambda$ with $n$ columns.
Denote by $T=\{(i,t_i)\}_{i=1}^{n} \in \mathcal{T}_n$ the transversal  in which the square $(i,t_i)$ is filled with a $1$ for all $i\leq n$. For example, the transversal $T=\{(1,6),  (2,2),  (3,8), (4, 4), (5,7), (6,1), (7,5), (8,3) \}$ of a Young diagram $\lambda=(8,8,8,8,8,5,5,5)$ is illustrated as   Figure \ref{transversal}.
	
	In this section, we will  consider permutations as permutation matrices. Given a permutation $\pi=\pi_1\pi_2\cdots\pi_n\in \mathcal{S}_n$, its corresponding {\em permutation matrix} is a transversal of the square Young diagram  $(\lambda_1, \lambda_2, \ldots, \lambda_n)$ with  $\lambda_1=\lambda_2=\cdots=\lambda_n=n$ in which the square   $(i, \pi_i)$ is filled with  a $1$ for all $1\leq i\leq n$ and all the other squares  are filled with $0's$.
	
	The   notion of pattern  avoidance is extended to transversals of a Young diagram in \cite{BW}  and \cite{BWX}.
	Given a permutation $\alpha$ of $\mathcal{S}_m$,  let   $M$ be its permutation matrix. A transversal $T$ of a Young diagram $\lambda$ with $n$ columns will be said to contain   $\alpha$ if there exists two subsets of the index set $[n]$, namely,  $R=\{r_1<r_2<\cdots<r_m\}$ and $C=\{c_1< c_2< \cdots< c_m\}$, such that  the   matrix restricted on rows $R$ and columns $C$ is a copy of $M$ and each of the squares $(c_i,r_j)$ falls within the Young diagram $\lambda$.   In this context,  the  permutation $\alpha$ is  called a pattern. For example, the transversal in  Figure \ref{transversal} contains the pattern $123$,   but avoids the pattern $321$. Denote by $\mathcal{T}_{\lambda}(\alpha)$ the set of all transversals of the  Young diagram $\lambda$ that avoid  the pattern $\alpha$.  Given two patterns $\alpha$ and $\beta$,   we say that $\alpha$ and $\beta$  are {\em shape-Wilf-equivalent}, denoted by $\alpha\sim_{s} \beta$,  if for any Young diagram $\lambda$ we have $|\mathcal{T}_{\lambda}(\alpha)|=|\mathcal{T}_{\lambda}(\beta)|$.

	    Backelin-West-Xin \cite{BWX} proved the following shape-Wilf-equivalence  for transversals of Young diagrams.

	\begin{theorem}\label{BWX2} (\cite{BWX}, Proposition  2.2)
		For   $k\geq 1$, we have $I_k\sim_{s} J_k$.
	\end{theorem}

	\begin{theorem}\label{BWX1} (\cite{BWX}, Proposition  2.3)
		For  any  patterns $\alpha$, $\beta$ and $\sigma$, if $\alpha\sim_{s}\beta$, then  $\alpha\oplus \sigma \sim_{s} \beta\oplus \sigma$.
	\end{theorem}
	
In order to prove Theorem \ref{generalth}, we shall derive   two theorems which are analogues  of Theorems \ref{BWX2} and  \ref{BWX1}. Before we state our results, we introduce some necessary definitions and notations.

	For  a Young diagram $\lambda$ with $k$ columns, we denote by $c_i(\lambda)$ the number of squares in column $i$ for all $i=1, 2, \ldots, k$.
	Given a transversal $T=\{(i,t_i)\}_{i=1}^{k}$ of the  Young diagram $\lambda$, an index $i$ $(2\leq i\leq k-1)$  is said to be a {\em peak} of $T$ if $c_{i-1}(\lambda)=c_{i}(\lambda)=c_{i+1}(\lambda)$ and $t_{i-1}<t_i>t_{i+1}$.  Denote by $\mathrm{Peak}(T)$ the set of peaks  of $T$.
	For example, if we let $T$ be the transversal indicated in Figure \ref{transversal}, then we have $\mathrm{Peak}(T)=\{3, 7\}$.
When restricted to permutation matrix, the peaks of transversals coincide with
the peaks of permutations.	
	
	For a self-conjugate Young diagram $\lambda$, a transversal $T$ of $\lambda$ is said to be   {\em symmetric} if the square $(i,j)$ is filled with a $1$ if and only if the square $(j,i)$ is filled with a $1$.  See Figure \ref{transversal} for an illustration of a symmetric transversal.  It is easily seen that the  permutation matrix of an involution in $\mathcal{I}_n$ is a symmetric  transversal of the square Young diagram $\lambda=(\lambda_1, \lambda_2, \ldots, \lambda_n)$ with $\lambda_1=\lambda_2=\cdots=\lambda_n=n$.  Denote by $\mathcal{ST}_{\lambda}(\alpha)$ the set of all symmetric transversals of the self-conjugate  Young diagram $\lambda$ that avoid the pattern $\alpha$.
	
Now we are ready to state our  two main theorems which will play essential  roles in the proof of
 Theorem \ref{generalth}.
	\begin{theorem}\label{mainth2}
		 For any self-conjugate Young diagram $\lambda$, there is a bijection $\Psi$ between $\mathcal{ST}_{\lambda}(J_3)$ and $\mathcal{ST}_{\lambda}(I_3)$ such that for any $T\in \mathcal{ST}_{\lambda}(J_3)$, we have $\mathrm{Peak}(T)=\mathrm{Peak}(\Psi(T))$.
	\end{theorem}

\begin{theorem}\label{mainth1}
		For  any  pattern   $\tau$ and any self-conjugate Young diagram $\lambda$, there exists  a bijection $\Phi$ between $\mathcal{ST}_{\lambda}(J_3\oplus \tau)$ and $\mathcal{ST}_{\lambda}(I_3\oplus \tau )$ such that for any $T\in \mathcal{ST}_{\lambda}(J_3\oplus \tau )$,  we have $\mathrm{Peak}(T)=\mathrm{Peak}(\Phi(T))$.
	\end{theorem}

\subsection{Proof of Theorem \ref{mainth2} }
	
\begin{figure}[h]
	\centering
	\begin{tikzpicture}[font =\small , scale = 0.4]
		\draw (-8,0)rectangle(-3.5,2);
		\draw [->] (-3,0.55)--(-0.5,0.55);
		\draw [<-] (-3,1.55)--(-0.5,1.55);
		\draw (0,0)rectangle(4.5,2);
		\draw [->] (5,0.55)--(7.5,0.55);
		\draw [<-] (5,1.55)--(7.5,1.55);
		\draw (8,0)rectangle(12.5,2);
		\draw [->] (13,0.55)--(15.5,0.55);
		\draw [<-] (13,1.55)--(15.5,1.55);
		\draw (16,0)rectangle(20.5,2);
		
		\draw [dashed] [->] (-6.25,2.5)--(-6.25,7.5);
		\draw [dashed] [<-] (-5.25,2.5)--(-5.25,7.5);
		\draw [<-] (17.75,2.5)--(17.75,7.5);
		\draw [->] (18.75,2.5)--(18.75,7.5);
		
		\draw (-8,8)rectangle(-3.5,10);
		\draw [<-] (-3,8.55)--(-0.5,8.55);
		\draw [->] (-3,9.55)--(-0.5,9.55);
		\draw (0,8)rectangle(4.5,10);
		\draw [<-] (5,8.55)--(7.5,8.55);
		\draw [->] (5,9.55)--(7.5,9.55);
		\draw (8,8)rectangle(12.5,10);
		\draw [<-] (13,8.55)--(15.5,8.55);
		\draw [->] (13,9.55)--(15.5,9.55);
		\draw (16,8)rectangle(20.5,10);

		\node at(-5.75,1){\footnotesize $\mathcal{ST}_\lambda(I_3)$};
		\node at(2.25,1){\footnotesize $\mathcal{SNM}_3(n)$};
		\node at(10.25,1){\footnotesize $\mathcal{SOC}_2(n)$};
		\node at(18.25,1){\footnotesize $\mathcal{SCN}(n)$};
		\node at(-5.75,9){\footnotesize $\mathcal{ST}_\lambda(J_3)$};
		\node at(2.25,9){\footnotesize $\mathcal{SCM}_3(n)$};
		\node at(10.25,9){\footnotesize $\mathcal{SOR}_2(n)$};
		\node at(18.25,9){\footnotesize $\mathcal{SCN}(n)$};
		
		\node at(-7.25,5){$\Psi^{-1}$};
		\node at(-4.5,5){$\Psi$};
		\node at(17,5){$\theta$};
		\node at(20,5){$\theta'$};
		\node at(-1.75,-0.25){\footnotesize $\chi$};
		\node at(-1.75,2.25){\footnotesize $\chi'$};
		\node at(6.25,-0.25){\footnotesize $\phi$};
		\node at(6.75,2.25){\footnotesize $\phi^{-1}$};
		\node at(14.55,-0.25){\footnotesize $\bar{\psi}$};
		\node at(14.25,2.25){\footnotesize $\bar{\psi}^{-1}$};
		\node at(-1.75,7.75){\footnotesize $\chi'$};
		\node at(-1.75,10.25){\footnotesize $\chi$};
		\node at(6.75,7.75){\footnotesize $\phi^{-1}$};
		\node at(6.25,10.25){\footnotesize $\phi$};
		\node at(14.55,7.75){\footnotesize ${\psi}^{-1}$};
		\node at(14.25,10.25){\footnotesize ${\psi}$};
								
	\end{tikzpicture}
\caption{A diagrammatic summary of the sets and bijections.} \label{fig:struct}
\end{figure}
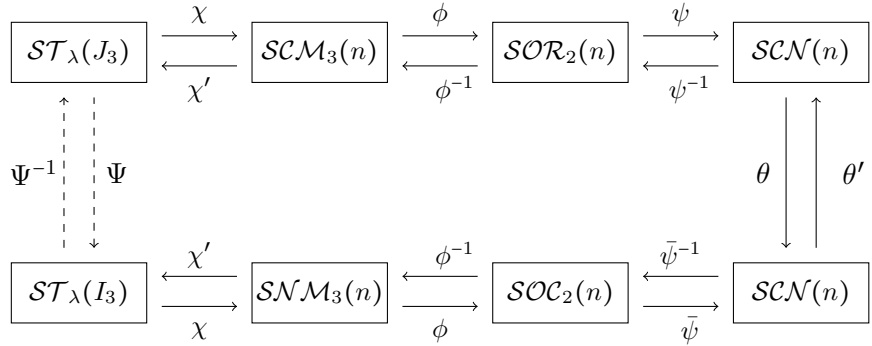

This subsection is devoted to the proof of Theorem \ref{mainth2}.  To this end,  we shall  construct a series of bijections and employ two bijections established by Chen-Deng-Du-Stanley-Yan \cite{ChenD}.  Figure \ref{fig:struct} outlines our sets
and maps of interest.

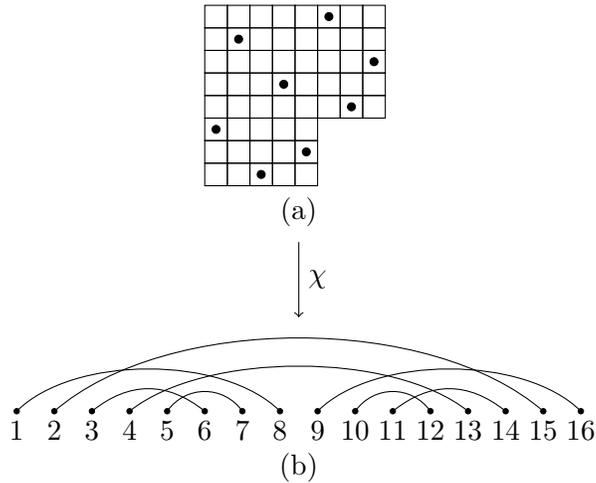
\begin{figure}[H]
\begin{center}
	\begin{tikzpicture}[font =\small , scale = 0.5]
		\draw (5,0)rectangle(5.6,0.6);
		\draw (5,0.6)rectangle(5.6,1.2);
		\draw (5,1.2)rectangle(5.6,1.8);
		\draw (5,1.8)rectangle(5.6,2.4);
		\draw (5,2.4)rectangle(5.6,3);
		\draw (5,3)rectangle(5.6,3.6);
		\draw (5,3.6)rectangle(5.6,4.2);
		\draw (5,4.2)rectangle(5.6,4.8);
		\draw (5.6,0)rectangle(6.2,0.6);
		\draw (5.6,0.6)rectangle(6.2,1.2);
		\draw (5.6,1.2)rectangle(6.2,1.8);
		\draw (5.6,1.8)rectangle(6.2,2.4);
		\draw (5.6,2.4)rectangle(6.2,3);
		\draw (5.6,3)rectangle(6.2,3.6);
		\draw (5.6,3.6)rectangle(6.2,4.2);
		\draw (5.6,4.2)rectangle(6.2,4.8);
		\draw (6.2,0)rectangle(6.8,0.6);
		\draw (6.2,0.6)rectangle(6.8,1.2);
		\draw (6.2,1.2)rectangle(6.8,1.8);
		\draw (6.2,1.8)rectangle(6.8,2.4);
		\draw (6.2,2.4)rectangle(6.8,3);
		\draw (6.2,3)rectangle(6.8,3.6);
		\draw (6.2,3.6)rectangle(6.8,4.2);
		\draw (6.2,4.2)rectangle(6.8,4.8);
    	\draw (6.8,0)rectangle(7.4,0.6);
		\draw (6.8,0.6)rectangle(7.4,1.2);
		\draw (6.8,1.2)rectangle(7.4,1.8);
		\draw (6.8,1.8)rectangle(7.4,2.4);
		\draw (6.8,2.4)rectangle(7.4,3);
		\draw (6.8,3)rectangle(7.4,3.6);
		\draw (6.8,3.6)rectangle(7.4,4.2);
		\draw (6.8,4.2)rectangle(7.4,4.8);
		\draw (7.4,0)rectangle(8,0.6);
		\draw (7.4,0.6)rectangle(8,1.2);
		\draw (7.4,1.2)rectangle(8,1.8);
 	    \draw (7.4,1.8)rectangle(8,2.4);
		\draw (7.4,2.4)rectangle(8,3);
		\draw (7.4,3)rectangle(8,3.6);
		\draw (7.4,3.6)rectangle(8,4.2);
		\draw (7.4,4.2)rectangle(8,4.8);
		\draw (8,1.8)rectangle(8.6,2.4);
		\draw (8,2.4)rectangle(8.6,3);
		\draw (8,3)rectangle(8.6,3.6);
		\draw (8,3.6)rectangle(8.6,4.2);
		\draw (8,4.2)rectangle(8.6,4.8);
		\draw (8.6,1.8)rectangle(9.2,2.4);
		\draw (8.6,2.4)rectangle(9.2,3);
		\draw (8.6,3)rectangle(9.2,3.6);
		\draw (8.6,3.6)rectangle(9.2,4.2);
		\draw (8.6,4.2)rectangle(9.2,4.8);
		\draw (9.2,1.8)rectangle(9.8,2.4);
		\draw (9.2,2.4)rectangle(9.8,3);
		\draw (9.2,3)rectangle(9.8,3.6);
		\draw (9.2,3.6)rectangle(9.8,4.2);
		\draw (9.2,4.2)rectangle(9.8,4.8);
		\node at(7.5,-0.7){(a)};
		\filldraw[black](5.3,1.5)circle(3pt);
		\filldraw[black](5.9,3.9)circle(3pt);
		\filldraw[black](6.5,0.3)circle(3pt);
		\filldraw[black](7.1,2.7)circle(3pt);
		\filldraw[black](7.7,0.9)circle(3pt);
		\filldraw[black](8.3,4.5)circle(3pt);
		\filldraw[black](8.9,2.1)circle(3pt);
		\filldraw[black](9.5,3.3)circle(3pt);
		
		\draw[->](7.5,-1.5)--(7.5,-3.5);
		\node at(8,-2.5){$\chi$};
		
		\filldraw[black](0,-6)circle(2pt);
		\node at(0,-6.5) {1};
		\filldraw[black](1,-6)circle(2pt);
		\node at(1,-6.5) {2};
		\filldraw[black](2,-6)circle(2pt);
		\node at(2,-6.5) {3};
		\filldraw[black](3,-6)circle(2pt);
		\node at(3,-6.5) {4};
		\filldraw[black](4,-6)circle(2pt);
		\node at(4,-6.5) {5};
		\filldraw[black](5,-6)circle(2pt);
		\node at(5,-6.5) {6};
		\filldraw[black](6,-6)circle(2pt);
		\node at(6,-6.5) {7};
		\filldraw[black](7,-6)circle(2pt);
		\node at(7,-6.5) {8};
		\filldraw[black](8,-6)circle(2pt);
		\node at(8,-6.5) {9};
		\filldraw[black](9,-6)circle(2pt);
		\node at(9,-6.5) {10};
		\filldraw[black](10,-6)circle(2pt);
		\node at(10,-6.5) {11};
		\filldraw[black](11,-6)circle(2pt);
		\node at(11,-6.5) {12};
		\filldraw[black](12,-6)circle(2pt);
		\node at(12,-6.5) {13};
		\filldraw[black](13,-6)circle(2pt);
		\node at(13,-6.5) {14};
		\filldraw[black](14,-6)circle(2pt);
		\node at(14,-6.5) {15};
		\filldraw[black](15,-6)circle(2pt);
		\node at(15,-6.5) {16};
		\draw (0,-6)..controls(1.5,-4.5)and(5.5,-4.5)..(7,-6);
		\draw (1,-6)..controls(4,-3.4)and(11,-3.4)..(14,-6);
		\draw (2,-6)..controls(3,-5.2)and(4,-5.2)..(5,-6);
		\draw (3,-6)..controls(5,-4.4)and(10,-4.4)..(12,-6);
		\draw (4,-6)..controls(4.5,-5.3)and(5.5,-5.3)..(6,-6);
		\draw (8,-6)..controls(9.5,-4.5)and(13.5,-4.5)..(15,-6);
		\draw (9,-6)..controls(9.5,-5.3)and(10.5,-5.3)..(11,-6);
		\draw (10,-6)..controls(11,-5.2)and(12,-5.2)..(13,-6);
		\node at(7.5,-7.5){(b)};
						
	\end{tikzpicture}
	
\end{center}
\caption{ An example of the bijection $\chi$ between $\mathcal{ST}_{\lambda}(J_3)$ and $\mathcal{SCM}_3(8)$ where $\lambda=(8,8,8,8,8,5,5,5)$.}\label{nc}
\end{figure}

Recall that a {\em set partition} $P$ of $[n]$ is a collection of nonempty subsets $\{B_1, B_2,$ $ \ldots, B_k\}$,  whose disjoint union is $[n]$. Each $B_i$ is called a {\em block} of $P$.
A (complete) {\em matching}  $M$ of $[2n] = \{1,2,\ldots, 2n\}$ is a  set partition of $[2n]$ into $n$ blocks and each block contains exactly two elements. Every set partition $P$ of $[n]$   can be represented
by a diagram with vertices drawn on a horizontal line in increasing order  and draw an arc connecting $i$ and  $j$ whenever $i$ and $j$ are   (numerically) consecutive in a block of $P$.  Such a diagram is called
the {\em linear representation} of $P$, see Figure \ref{nc}(b) for an example.  Throughout the paper, an arc  is always written as $(i,j)$ with $i<j$. For an arc $(i,j)$,  we call $i$ and $j$  its {\em left-hand endpoint} and {\em right-hand endpoint}, respectively.  A  vertex $i$ is said to an {\em opener} of a matching $M$ if $i$ is a left-hand endpoint of some arc. Similarly,  a vertex $i$ is said to a {\em closer} of a matching $M$ if $i$ is a right-hand endpoint of some arc.
For a set partition $P$,  we say that $k$ arcs $(i_1, j_1), (i_2, j_2), \ldots, (i_k, j_k)$ form a
$k$-crossing if $i_1<i_2<\cdots <i_k<j_1<j_2<\cdots<j_k$, and let $\mathrm{cr}(P)$ denote the maximal $k$ such that $P$ has a $k$-crossing.  A set partition without any $k$-crossings  is said to be {\em $k$-noncrossing}.
Similarly,  a {$k$-nesting} is a set of $k$ arcs $(i_1, j_1), (i_2, j_2), \ldots, (i_k, j_k)$   such that  $i_1<i_2<\cdots <i_k<j_k<\cdots<j_2<j_1$, and  denote by $\mathrm{ne}(P)$   the maximal $k$ such that $P$ has a $k$-nesting.
 A set partition without any $k$-nestings  is said to be {\em $k$-nonnesting}. Chen-Deng-Du-Stanley-Yan \cite{ChenD} proved that $k$-nonnesting set partitions of $[ n]$ are equinumerous with $k$-noncrossing set partitions  of $[n]$ bijectively using vacillating tableaux as an intermediate object.

For a set partition $P$ of $[n]$, let $P^{r}$ denote the set partition obtained from $P$ by reflecting its linear representation in the vertical line $x={n+1\over 2}$. Equivalently, $(i,j)$ is an arc of $P$ if and only if $(n+1-j, n+1-i)$ is an arc of $P^{r}$. A set partition $P$ is said to be {\em bilaterally symmetric} if $P^{r}=P$.
A  {\em vacillating   tableau} of shape $\lambda$ and length $2n$ is a sequence $(\emptyset=\lambda^0 , \lambda^1, \ldots, \lambda^{2n}=\lambda)$ of integer partitions such that  (\rmnum{1})  $\lambda^{2i+1}$ is obtained from $\lambda^{2i}$ by doing nothing or deleting a square, (\rmnum{2})  $\lambda^{2i}$ is obtained from $\lambda^{2i-1}$ by adding a square or doing nothing.    An {\em oscillating    tableau} of shape $\lambda$ and length $n$ is a sequence $(\emptyset=\lambda^0, \lambda^1, \ldots, \lambda^{n}=\lambda)$ of integer partitions such that     $\lambda^{i}$ is obtained from $\lambda^{i-1}$ by   adding a square or deleting a square.
For example, $(\emptyset, (1), (1,1), (2,1), (1,1), (2,1), (2), (1), \emptyset)$  is an oscillating    tableau of length $8$ and  shape $\emptyset$.
In what follows, vacillating   tableaux and oscillating    tableaux are always of shape $\emptyset$ unless specified otherwise.

For a transversal  $T$ of
  a Young diagram $\lambda$, let $\mathrm{type}(T)$ denote  the sequence  obtained from  $T$ by  tracing the right  boundary of $\lambda$ from  bottom-left to top-right and writing  $U$ (resp. $D$) whenever we encounter a east step (resp. a north step). For example, let $T$ be the transversal    as illustrated in Figure \ref{nc}(a), we have  $\mathrm{type(T)}=UUUUUDDDUUUDDDDD$.
For a matching $M$, we let $\mathrm{type}(M)$ denote the  sequence  obtained from $M$ by tracing the vertices of matching from left to right and writing  $U$ (resp. $D$) whenever we encounter an opener (resp. a closer). For example, let $M$ be the matching  as illustrated in Figure \ref{nc}(b), then we have $\mathrm{type(M)}=UUUUUDDDUUUDDDDD$.

%Let $\mathcal{T}_n$ denote the set of transversals of Young diagrams with $n$ columns and
Let $\mathcal{M}_n$  denote the set of matchings of $[2n]$.
There is   a simple  map $\chi: \mathcal{T}_n\rightarrow \mathcal{M}_n$.    Given a transversal $T \in \mathcal{T}_n$, we construct a matching $M=\chi(T)$ satisfying that
\begin{itemize}
\item   $\mathrm{type}(M)=\mathrm{type}(T)$;
\item   there is an arc connecting the $i$-th left-to-right opener and the $j$-th right-to-left closer whenever   there is a $1$ positioned at the square $(i, j)$.
\end{itemize}

Conversely, given a matching $M\in \mathcal{M}_n$, we recover a transversal $T = \chi'(M)$ as follows:
 \begin{itemize}
\item  $\mathrm{type}(T)=\mathrm{type}(M)$;
\item  Put a $1$ in the square $(i, j)$ whenever there is an arc connecting the $i$-th left-to-right opener and the $j$-th righ-to-left closer.
\end{itemize}
Obviously, the maps $\chi$ and $\chi'$ are inverses of each other, and hence the map $\chi$ is a bijection. It is routine to check that a symmetric transversal corresponds to a bilaterally  symmetric matching.   Moreover, a $k$-crossing (resp.  $k$-nesting) in a matching corresponds to a  pattern $J_k$ (resp. $I_k$) in the corresponding transversal.

  In a matching $M$, three arcs $(i_1, j_1)$,  $(i_2, j_2)$ and $(i_3,j_3)$ are said to form a {\em valley} of   $M$ if $i_2=i_1+1$, $i_3=i_2+1$,  and $j_1>j_2<j_3$.   In this context,    $i_2$     is said to be a {\em valley index}   of $M$.     Denote by $\mathrm{Val}(M)$ the set of valley indexes  of $M$. For a matching $M$, assume that $\mathrm{Val}(M)=\{a_1, a_2, \ldots, a_k\}$. If $a_i$ is the $b_i$-th left-to-right opener of $M$, define
$\widetilde{\mathrm{Val }}(M)=\{b_1, b_2, \ldots, b_k\}$.  For example, if we let $M$ be the matching as shown in Figure \ref{nc}(b), then we have $\mathrm{Val}(M)=\{3,10\}$ and $\widetilde{\mathrm{Val }}(M)=\{3,7\}$.
% It is apparent that  for a matching $M$,  $\widetilde{\mathrm{Val }}(M)$ is uniquely determined by $\mathrm{type}(M)$ and $\mathrm{Val}(M)$.
It is apparent that  for a matching $M$ of a given type,
$\widetilde{\mathrm{Val }}(M)$ is uniquely determined by $\mathrm{Val}(M)$, and vise versa.

 According to the construction of $\chi$, one can easily  verify that $i\in \mathrm{Peak}(T)$ if and only if $i\in \widetilde{ \mathrm{Val}}(\chi(T))$. Let $\mathcal{SCM}_k(n)$ and $\mathcal{SNM}_k(n)$ denote the sets of $k$-noncrossing  and  $k$-nonnesting bilaterally  symmetric matchings of $[2n]$, respectively.
 The following two properties of $\chi$ can be deduced by the above analysis.

 \begin{lemma}\label{lemchi1}
Let $\lambda$ be a Young diagram with $n$ columns.
 The map $\chi$ induces a bijection between $\mathcal{ST}_{\lambda}(J_k)$ and $\mathcal{SCM}_k(n)$ such that for any transversal $T\in \mathcal{ST}_{\lambda}(J_k)$, we have $\mathrm{type}(T)=\mathrm{type}(\chi(T))$ and $ \mathrm{Peak}(T)=\widetilde{\mathrm{Val}}(\chi(T))$.
 \end{lemma}

\begin{lemma}\label{lemchi2}

Let $\lambda$ be a Young diagram with $n$ columns.
 The map $\chi$ induces a bijection between $\mathcal{ST}_{\lambda}(I_k)$ and $\mathcal{SNM}_k(n)$ such that for any transversal $T\in \mathcal{ST}_{\lambda}(I_k)$, we have $\mathrm{type}(T)=\mathrm{type}(\chi(T))$ and $ \mathrm{Peak}(T)=\widetilde{\mathrm{Val}}(\chi(T))$.
 \end{lemma}

For example, let $T\in \mathcal{ST}_{\lambda}(J_3)$ be a  transversal as illustrated in Figure \ref{nc}(a), where $\lambda=(8,8,8,8,8,5,5,5)$.  By applying the bijection $\chi$ to $T$, we get a matching $\chi(T)\in  \mathcal{SCM}_3(8)$  as shown in Figure \ref{nc}(b).  One can easily check that $\mathrm{type}(T)=\mathrm{type}(\chi(T))=UUUUUDDDUUUDDDDD$ and
$ \mathrm{Peak}(T)=\widetilde{\mathrm{Val}}(\chi(T))=\{3,7\}$.

\begin{figure}[H]
\begin{center}
	\begin{tikzpicture}[font =\small , scale = 0.5]
  		\draw (5,0)rectangle(5.6,0.6);
  		\draw (5,0.6)rectangle(5.6,1.2);
  		\draw (5,1.2)rectangle(5.6,1.8);
  		\draw (5,1.8)rectangle(5.6,2.4);
  		\draw (5,2.4)rectangle(5.6,3);
  		\draw (5,3)rectangle(5.6,3.6);
  		\draw (5,3.6)rectangle(5.6,4.2);
  		\draw (5,4.2)rectangle(5.6,4.8);
  		\draw (5.6,0)rectangle(6.2,0.6);
  		\draw (5.6,0.6)rectangle(6.2,1.2);
  		\draw (5.6,1.2)rectangle(6.2,1.8);
  		\draw (5.6,1.8)rectangle(6.2,2.4);
  		\draw (5.6,2.4)rectangle(6.2,3);
  		\draw (5.6,3)rectangle(6.2,3.6);
  		\draw (5.6,3.6)rectangle(6.2,4.2);
  		\draw (5.6,4.2)rectangle(6.2,4.8);
  		\draw (6.2,0)rectangle(6.8,0.6);
  		\draw (6.2,0.6)rectangle(6.8,1.2);
  		\draw (6.2,1.2)rectangle(6.8,1.8);
  		\draw (6.2,1.8)rectangle(6.8,2.4);
  		\draw (6.2,2.4)rectangle(6.8,3);
  		\draw (6.2,3)rectangle(6.8,3.6);
  		\draw (6.2,3.6)rectangle(6.8,4.2);
  		\draw (6.2,4.2)rectangle(6.8,4.8);
  		\draw (6.8,0)rectangle(7.4,0.6);
  		\draw (6.8,0.6)rectangle(7.4,1.2);
  		\draw (6.8,1.2)rectangle(7.4,1.8);
  		\draw (6.8,1.8)rectangle(7.4,2.4);
  		\draw (6.8,2.4)rectangle(7.4,3);
  		\draw (6.8,3)rectangle(7.4,3.6);
  		\draw (6.8,3.6)rectangle(7.4,4.2);
  		\draw (6.8,4.2)rectangle(7.4,4.8);
  		\draw (7.4,0)rectangle(8,0.6);
  		\draw (7.4,0.6)rectangle(8,1.2);
  		\draw (7.4,1.2)rectangle(8,1.8);
  		\draw (7.4,1.8)rectangle(8,2.4);
  		\draw (7.4,2.4)rectangle(8,3);
  		\draw (7.4,3)rectangle(8,3.6);
  		\draw (7.4,3.6)rectangle(8,4.2);
  		\draw (7.4,4.2)rectangle(8,4.8);
  		\draw (8,1.8)rectangle(8.6,2.4);
  		\draw (8,2.4)rectangle(8.6,3);
  		\draw (8,3)rectangle(8.6,3.6);
  		\draw (8,3.6)rectangle(8.6,4.2);
  		\draw (8,4.2)rectangle(8.6,4.8);
  		\draw (8.6,1.8)rectangle(9.2,2.4);
  		\draw (8.6,2.4)rectangle(9.2,3);
  		\draw (8.6,3)rectangle(9.2,3.6);
  		\draw (8.6,3.6)rectangle(9.2,4.2);
  		\draw (8.6,4.2)rectangle(9.2,4.8);
  		\draw (9.2,1.8)rectangle(9.8,2.4);
  		\draw (9.2,2.4)rectangle(9.8,3);
  		\draw (9.2,3)rectangle(9.8,3.6);
  		\draw (9.2,3.6)rectangle(9.8,4.2);
  		\draw (9.2,4.2)rectangle(9.8,4.8);
  		\node at(7.5,-0.7){(a)};
  		\filldraw[black](5.3,1.5)circle(3pt);
  		\filldraw[black](5.9,2.7)circle(3pt);
  		\filldraw[black](6.5,0.3)circle(3pt);
  		\filldraw[black](7.1,3.9)circle(3pt);
  		\filldraw[black](7.7,0.9)circle(3pt);
  		\filldraw[black](8.3,4.5)circle(3pt);
  		\filldraw[black](8.9,2.1)circle(3pt);
  		\filldraw[black](9.5,3.3)circle(3pt);
  		
  		\draw[->](7.5,-1.5)--(7.5,-3.5);
  		\node at(8,-2.5){$\chi$};
  		
  		\filldraw[black](0,-6)circle(2pt);
  		\node at(0,-6.5) {1};
  		\filldraw[black](1,-6)circle(2pt);
  		\node at(1,-6.5) {2};
  		\filldraw[black](2,-6)circle(2pt);
  		\node at(2,-6.5) {3};
  		\filldraw[black](3,-6)circle(2pt);
  		\node at(3,-6.5) {4};
  		\filldraw[black](4,-6)circle(2pt);
  		\node at(4,-6.5) {5};
  		\filldraw[black](5,-6)circle(2pt);
  		\node at(5,-6.5) {6};
  		\filldraw[black](6,-6)circle(2pt);
  		\node at(6,-6.5) {7};
  		\filldraw[black](7,-6)circle(2pt);
  		\node at(7,-6.5) {8};
  		\filldraw[black](8,-6)circle(2pt);
  		\node at(8,-6.5) {9};
  		\filldraw[black](9,-6)circle(2pt);
  		\node at(9,-6.5) {10};
  		\filldraw[black](10,-6)circle(2pt);
  		\node at(10,-6.5) {11};
  		\filldraw[black](11,-6)circle(2pt);
  		\node at(11,-6.5) {12};
  		\filldraw[black](12,-6)circle(2pt);
  		\node at(12,-6.5) {13};
  		\filldraw[black](13,-6)circle(2pt);
  		\node at(13,-6.5) {14};
  		\filldraw[black](14,-6)circle(2pt);
  		\node at(14,-6.5) {15};
  		\filldraw[black](15,-6)circle(2pt);
  		\node at(15,-6.5) {16};
  		\draw (0,-6)..controls(1.5,-4.5)and(5.5,-4.5)..(7,-6);
  		\draw (1,-6)..controls(3,-3.4)and(10,-3.4)..(12,-6);
  		\draw (2,-6)..controls(3,-5.2)and(4,-5.2)..(5,-6);
  		\draw (3,-6)..controls(5,-3.4)and(12,-3.4)..(14,-6);
  		\draw (4,-6)..controls(4.5,-5.3)and(5.5,-5.3)..(6,-6);
  		\draw (8,-6)..controls(9.5,-4.5)and(13.5,-4.5)..(15,-6);
  		\draw (9,-6)..controls(9.5,-5.3)and(10.5,-5.3)..(11,-6);
  		\draw (10,-6)..controls(11,-5.2)and(12,-5.2)..(13,-6);
  		\node at(7.5,-7.5){(b)};
  		
  \end{tikzpicture}
	
\end{center}
\caption{ An example of the bijection $\chi$ between $\mathcal{ST}_{\lambda}(I_3)$ and $\mathcal{SNM}_3(8)$ where $\lambda=(8,8,8,8,8,5,5,5)$.}\label{nn}
\end{figure}
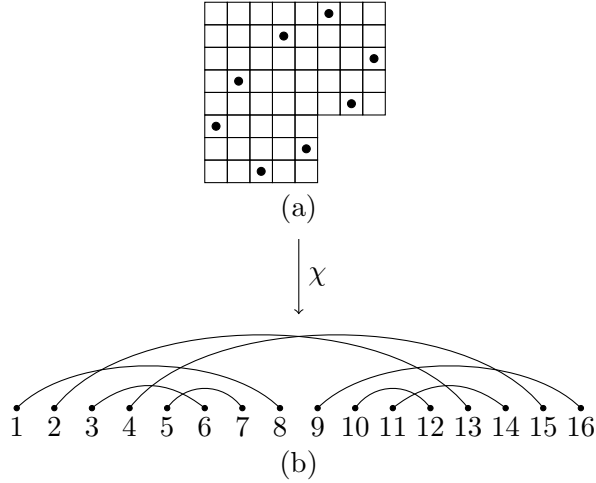

For example, let $T\in \mathcal{ST}_{\lambda}(I_3)$ be a  transversal as illustrated in Figure \ref{nn}(a), where $\lambda=(8,8,8,8,8,5,5,5)$.  By applying the bijection $\chi$ to $T$, we get a matching $\chi(T)\in  \mathcal{SNM}_3(8)$  as shown in Figure \ref{nn}(b).  One can easily check that $\mathrm{type}(T)=\mathrm{type}(\chi(T))=UUUUUDDDUUUDDDDD$ and
$ \mathrm{Peak}(T)=\widetilde{\mathrm{Val}}(\chi(T))=\{3,7\}$.

In \cite{ChenD}, Chen-Deng-Du-Stanley-Yan established a bijection $\phi$ between set partitions   and  vacillating   tableaux. When  restricting to matchings, the bijection $\phi$ induces a bijection between matchings and  oscillating    tableaux.  The explicit description of the bijection between matchings and  oscillating    tableaux can also be found in \cite{StanleyVol2} and \cite{Sundaram}.
Chen-Deng-Du-Stanley-Yan \cite{ChenD} proved that the bijection $\phi$ verifies  the following celebrated property.
 \begin{theorem}\label{chenth6}
 (\cite{ChenD},  Theorem 6)
   Let $P$ be a set partition of $[n]$ and $\phi(P)=(\lambda^0, \lambda^1, \ldots, \lambda^{n})$.  Then  $\mathrm{cr}(P)$ is the most number
of rows in any $\lambda^i$, and $\mathrm{ne}(P)$ is the most number of columns in any $\lambda^i$.
 \end{theorem}

For a vacillating  (resp. an oscillating) tableau $V$, reading $V$ backward still gives a   vacillating (resp. an oscillating ) tableau, denoted by $V^{r}$. If $V^{r}=V$, we say that $V$ is {\em symmetric}.
By using Sch\"{u}tzenberger's theorem (see \cite{StanleyVol2}, Chapter 7.11) for the ordinary  RSK correspondence (see \cite{StanleyVol2}, Chapter 7.13), Xin-Zhang \cite{xinzhang} proved that the bijection $\phi$ has the following celebrated property.

\begin{theorem}\label{xinzhang}
(\cite{xinzhang},  Theorem 1)
For any given set partition $P$ and vacillating tableau $V$, $\phi(P^{r})=V^{r}$ if and only if $\phi(P)=V$.
\end{theorem}

In the following, we give a review of the bijection $\phi$ between  matchings of $[2n]$ and  oscillating  tableaux of length $2n$.

\noindent {\bf The bijection $\phi$ from Matchings to Oscillating Tableaux.}\\
 Given a matching $M$ of $[2n]$ with the
linear representation, we construct the sequence of SYT's, hence the oscillating tableau $\phi(M)$ as follows: Start from the empty SYT by letting $T_{2n}=\emptyset$, read the number $j\in [2n]$ one by one from
$2n$ to 1, and let $T_{j-1}$ be the SYT obtained from $T_j$ for each $j$ by the following procedure.
 \begin{itemize}
\item If  $j$ is the right-hand endpoint of an arc $(i,j)$, then insert $i$ (by the RSK algorithm) into the tableau.
\item If $j$ is the left-hand endpoint of an arc $(j, k)$,  then  remove $j$.
\end{itemize}
 Then the oscillating tableau $\phi(M)$  is the sequence  of shapes of the above SYT's.

For example, let $ M=\{(1, 8), (2,15), (3,6), (4, 13), (5,7), (9,16), (10,12), (11,14) \} $ be a $3$-noncrossing matching whose linear representation is   illustrated in Figure \ref{nc}(b). By applying the bijection $\phi$, its corresponding  oscillating tableau $\phi(M)$ is indicated in Table \ref{phi1}.

\begin{table}[!h]
\centering
\small
\renewcommand\arraystretch{0.75}
\caption{An oscillating tableau $\phi(M)$ corresponding to the matching  $M=\{(1, 8), (2,15), (3,6), (4, 13), (5,7), (9,16), (10,12), (11,14) \}$.}\label{phi1}
 \vskip 2mm

\begin{tabular}{|c|m{2cm}<{\centering}|c|c| m{2cm}<{\centering}|c|}
\hline
 $i$ &  \vskip 0.5mm  $T_i$       & $\lambda^i$      &  $i$   & \vskip 0.5mm $T_i$     & $\lambda^i$    \\ \hline
$0$ &  $\emptyset$   & $\emptyset$       & $9$    &
\vskip 1mm
$\begin{array}[b]{*{2}c}\cline{1-2}
\lr{2}&\lr{4}\\\cline{1-2}
\lr{9}\\\cline{1-1}
\end{array}$
& $(2,1)$                                    \\ \hline

 $1$  &
$\begin{array}[b]{*{1}c}\cline{1-1}
\lr{1}\\\cline{1-1}
\end{array}$
& $(1)$ & $10$ &  \vskip 1mm
$\begin{array}[b]{*{3}c}\cline{1-3}
\lr{2}&\lr{4}&\lr{10}\\\cline{1-3}
\lr{9}\\\cline{1-1}
\end{array}$
& $(3,1)$                                 \\ \hline

$2$  & \vskip 1mm
$\begin{array}[b]{*{1}c}\cline{1-1}
\lr{1}\\\cline{1-1}
\lr{2}\\\cline{1-1}
\end{array}$
& $(1,1)$  & $11$ &  \vskip 1mm
$\begin{array}[b]{*{3}c}\cline{1-3}
\lr{2}&\lr{4}&\lr{10}\\\cline{1-3}
\lr{9}&\lr{11}\\\cline{1-2}
\end{array}$
&   $(3,2)$                             \\ \hline

$3$  & \vskip 1mm
$\begin{array}[b]{*{2}c}\cline{1-2}
\lr{1}&\lr{3}\\\cline{1-2}
\lr{2}\\\cline{1-1}
\end{array}$
& $(2,1)$  & $12$  & \vskip 1mm
$\begin{array}[b]{*{2}c}\cline{1-2}
\lr{2}&\lr{4}\\\cline{1-2}
\lr{9}&\lr{11}\\\cline{1-2}
\end{array}$
&    $(2,2)$                            \\ \hline

$4$  & \vskip 1mm
$\begin{array}[b]{*{2}c}\cline{1-2}
\lr{1}&\lr{3}\\\cline{1-2}
\lr{2}&\lr{4}\\\cline{1-2}
\end{array}$
& $(2,2)$  & $13$  & \vskip 1mm
$\begin{array}[b]{*{2}c}\cline{1-2}
\lr{2}&\lr{11}\\\cline{1-2}
\lr{9}\\\cline{1-1}
\end{array}$
&    $(2,1)$                            \\ \hline

$5$  & \vskip 1mm
$\begin{array}[b]{*{3}c}\cline{1-3}
\lr{1}&\lr{3}&\lr{5}\\\cline{1-3}
\lr{2}&\lr{4}\\\cline{1-2}
\end{array}$
& $(3,2)$  & $14$  & \vskip 1mm
$\begin{array}[b]{*{2}c}\cline{1-1}
\lr{2}\\\cline{1-1}
\lr{9}\\\cline{1-1}
\end{array}$
&    $(1,1)$                            \\ \hline

$6$  & \vskip 1mm
$\begin{array}[b]{*{3}c}\cline{1-3}
\lr{1}&\lr{4}&\lr{5}\\\cline{1-3}
\lr{2}\\\cline{1-1}
\end{array}$
& $(3,1)$  & $15$  & \vskip 1mm
$\begin{array}[b]{*{2}c}\cline{1-1}
\lr{9}\\\cline{1-1}
\end{array}$
&    $(1)$                            \\ \hline

$7$  & \vskip 1mm
$\begin{array}[b]{*{2}c}\cline{1-2}
\lr{1}&\lr{4}\\\cline{1-2}
\lr{2}\\\cline{1-1}
\end{array}$
& $(2,1)$  & $16$  &  $\emptyset$ &  $\emptyset$        \\ \hline

$8$  & \vskip 1mm
$\begin{array}[b]{*{2}c}\cline{1-2}
\lr{2}&\lr{4}\\\cline{1-2}
\end{array}$
& $(2)$  &   &   &         \\ \hline
\end{tabular}
\end{table}

For example, let $ M=\{(1, 8), (2,13), (3,6), (4, 15), (5,7), (9,16), (10,12), (11,14) \} $ be a $3$-nonnesting matching whose linear representation is   illustrated in Figure \ref{nn}(b). By applying the bijection $\phi$, its corresponding  oscillating tableau $\phi(M)$ is indicated in Table \ref{phi2}.

\begin{table}[!h]
\centering
\small
\renewcommand\arraystretch{0.75}
\caption{An oscillating tableau $\phi(M)$ corresponding to the matching  $M=\{(1, 8), (2,13), (3,6), (4, 15), (5,7), (9,16), (10,12), (11,14) \}$.}\label{phi2}
 \vskip 2mm

\begin{tabular}{|c|m{2cm}<{\centering}|c|c| m{2cm}<{\centering}|c|}
\hline
$i$ &  \vskip 0.5mm  $T_i$       & $\lambda^i$      &  $i$   & \vskip 0.5mm $T_i$     & $\lambda^i$    \\ \hline
$0$ &  $\emptyset$   & $\emptyset$       & $9$    &
\vskip 1mm
$\begin{array}[b]{*{1}c}\cline{1-1}
\lr{2}\\\cline{1-1}
\lr{4}\\\cline{1-1}
\lr{9}\\\cline{1-1}
\end{array}$
& $(1,1,1)$                                    \\ \hline

 $1$  &
$\begin{array}[b]{*{1}c}\cline{1-1}
\lr{1}\\\cline{1-1}
\end{array}$
& $(1)$ & $10$ &  \vskip 1mm
$\begin{array}[b]{*{2}c}\cline{1-2}
\lr{2}&\lr{10}\\\cline{1-2}
\lr{4}\\\cline{1-1}
\lr{9}\\\cline{1-1}
\end{array}$
& $(2,1,1)$                                 \\ \hline

$2$  & \vskip 1mm
$\begin{array}[b]{*{1}c}\cline{1-1}
\lr{1}\\\cline{1-1}
\lr{2}\\\cline{1-1}
\end{array}$
& $(1,1)$  & $11$ &  \vskip 1mm
$\begin{array}[b]{*{2}c}\cline{1-2}
\lr{2}&\lr{10}\\\cline{1-2}
\lr{4}&\lr{11}\\\cline{1-2}
\lr{9} \\\cline{1-1}
\end{array}$
&   $(2,2,1)$                             \\ \hline

$3$  & \vskip 1mm
$\begin{array}[b]{*{2}c}\cline{1-2}
\lr{1}&\lr{3}\\\cline{1-2}
\lr{2}\\\cline{1-1}
\end{array}$
& $(2,1)$  & $12$  & \vskip 1mm
$\begin{array}[b]{*{2}c}\cline{1-2}
\lr{2}&\lr{11}\\\cline{1-2}
\lr{4}\\\cline{1-1}
\lr{9}\\\cline{1-1}
\end{array}$
&    $(2,1,1)$                            \\ \hline

$4$  & \vskip 1mm
$\begin{array}[b]{*{2}c}\cline{1-2}
\lr{1}&\lr{3}\\\cline{1-2}
\lr{2}\\\cline{1-1}
\lr{4}\\\cline{1-1}
\end{array}$
& $(2,1,1)$  & $13$  & \vskip 1mm
$\begin{array}[b]{*{2}c}\cline{1-2}
\lr{4}&\lr{11}\\\cline{1-2}
\lr{9}\\\cline{1-1}
\end{array}$
&    $(2,1)$                            \\ \hline

$5$  & \vskip 1mm
$\begin{array}[b]{*{2}c}\cline{1-2}
\lr{1}&\lr{3}\\\cline{1-2}
\lr{2}&\lr{5}\\\cline{1-2}
\lr{4}\\\cline{1-1}
\end{array}$
& $(2,2,1)$  & $14$  & \vskip 1mm
$\begin{array}[b]{*{2}c}\cline{1-1}
\lr{4}\\\cline{1-1}
\lr{9}\\\cline{1-1}
\end{array}$
&    $(1,1)$                            \\ \hline

$6$  & \vskip 1mm
$\begin{array}[b]{*{2}c}\cline{1-2}
\lr{1}&\lr{5}\\\cline{1-2}
\lr{2}\\\cline{1-1}
\lr{4}\\\cline{1-1}
\end{array}$
& $(2,1,1)$  & $15$  & \vskip 1mm
$\begin{array}[b]{*{2}c}\cline{1-1}
\lr{9}\\\cline{1-1}
\end{array}$
&    $(1)$                            \\ \hline

$7$  & \vskip 1mm
$\begin{array}[b]{*{1}c}\cline{1-1}
\lr{1}\\\cline{1-1}
\lr{2}\\\cline{1-1}
\lr{4}\\\cline{1-1}
\end{array}$
& $(1,1,1)$  & $16$  &  $\emptyset$ &  $\emptyset$        \\ \hline

$8$  & \vskip 1mm
$\begin{array}[b]{*{1}c}\cline{1-1}
\lr{2}\\\cline{1-1}
\lr{4}\\\cline{1-1}
\end{array}$
& $(1,1)$  &   &   &         \\ \hline
\end{tabular}
\end{table}

Let $\mathcal{OR}_{k}(n)$ denote the set of  oscillating  tableaux $(\lambda^0, \lambda^1, \ldots, $ $\lambda^{2n})$ in which  any $\lambda^i$ has at most $k$ rows.  Similarly, let $\mathcal{OC}_{k}(n)$ denote the set of oscillating  tableaux $(\lambda^0, \lambda^1, \ldots, \lambda^{2n})$ in which  any $\lambda^i$ has at most $k$ columns.
Let $\mathcal{SOR}_{k}(n)$ denote the set of  symmetric oscillating  tableaux of $\mathcal{OR}_{k}(n)$
and let $\mathcal{SOC}_{k}(n)$ denote the set of  symmetric oscillating  tableaux of $\mathcal{OC}_{k}(n)$.
%Let $\mathcal{SOR}_{k}(n)$ denote the set of  symmetric oscillating  tableaux $(\emptyset=\lambda^0, \lambda^1, \ldots, $ $\lambda^{2n}=\emptyset)$ in which  any $\lambda^i$ has at most $k$ rows.  Similarly, let $\mathcal{SOC}_{k}(n)$ denote the set of  symmetric oscillating  tableaux $(\emptyset=\lambda^0, \lambda^1, \ldots, \lambda^{2n}=\emptyset)$ in which  any $\lambda^i$ has at most $k$ columns.
For an oscillating tableau $O=(\lambda^0, \lambda^1, \ldots, \lambda^{2n})$, let $\mathrm{type}(O)$ denote the sequence  obtained from  $O$ by   reading $O$  forward and   writing   $U$  (resp.   $D$) whenever  $\lambda^i$ is obtained from $\lambda^{i-1}$ by adding (resp. deleting) a square.
For example, the type of the  oscillating  tableau in Table \ref{phi1} is given by $UUUUUDDDUUUDDDDD$.

For an  oscillating tableau $O=(\lambda^0, \lambda^1, \ldots, \lambda^{2n})\in \mathcal{SOR}_{2}(n)$,  an index $i$ $(1< i<2n)$    is said to be a {\em valley} of $O$ if (\rmnum{1}) $\lambda^{i-1}$ is obtained from $\lambda^{i-2}$ by adding  a square,  (\rmnum{2})   $\lambda^{i}$ is obtained from $\lambda^{i-1}$ by adding  a square at the first row,  and (\rmnum{3})   $\lambda^{i+1}$ is obtained from $\lambda^i$ by adding a square at the second row.
Denote by $\mathrm{Val}(O)$ the set of all valleys of $O$. For example, if we let $O\in \mathcal{SOR}_{2}(8)$ be the oscillating tableau as shown in Table  \ref{phi1},  we have $\mathrm{Val}(O)=\{3,10\}$.

For an  oscillating tableau $O=(\lambda^0, \lambda^1, \ldots, \lambda^{2n})\in \mathcal{SOC}_{2}(n)$,  an index $i$ $(1 < i<2n)$    is said to be a {\em peak} of $O$ if (\rmnum{1}) $\lambda^{i-1}$ is obtained from $\lambda^{i-2}$ by adding  a square at the first column,  (\rmnum{2})   $\lambda^{i}$ is obtained from $\lambda^{i-1}$ by adding  a square at the second column,  and (\rmnum{3})   $\lambda^{i+1}$ is obtained from $\lambda^i$ by adding a square.
Denote by $\mathrm{Peak}(O)$ the set of all peaks of $O$.  For example, if we let $O\in \mathcal{SOC}_{2}(8)$ be the oscillating tableau as shown in Table  \ref{phi2},  we have $\mathrm{Peak}(O)=\{3,10\}$.

\begin{lemma}\label{lemphi1}
The bijection $\phi$ induces a bijection between     $\mathcal{SCM}_3(n)$ and $\mathcal{SOR}_2(n)$ such that for any  matching $M\in \mathcal{SCM}_3(n)$, we have $\mathrm{type}(M)=\mathrm{type}(\phi(M))$ and $ \mathrm{Val}(M)= \mathrm{Val}(\phi(M))$.
\end{lemma}
\pf By Theorems \ref{chenth6} and \ref{xinzhang},  the bijection $\phi$ induces a bijection between $\mathcal{SCM}_3(n)$ and $\mathcal{SOR}_2(n)$.
  Given a matching $M\in \mathcal{SCM}_3(n)$, let $\phi(M)= (\lambda^0, \lambda^1, \ldots, \lambda^{2n})$.
  According to the construction of $\phi$, one can easily check that $\lambda^i$ is obtained from $\lambda^{i-1}$ by adding a square if and only if $i$ is an opener, and $\lambda^i$ is obtained from $\lambda^{i-1}$ by deleting a square if and only if $i$ is a closer. This ensures that $\mathrm{type}(M)=\mathrm{type}(\phi(M))$.

  Next we aim to show that $ \mathrm{Val}(M)= \mathrm{Val}(\phi(M))$.  By applying the map $\phi$ to $M$, we obtain a sequence $(T_0, T_1, \ldots, T_{2n})$ of SYT's  such that $T_i$ is of shape $\lambda^i$.  Assume that $k\in \mathrm{Val}(M)$. Then three arcs $(k-1, j_1), (k, j_2)$ and $(k+1, j_3)$ with $j_1>j_2<j_3$ form a valley of $M$. As $M$ is bilaterally symmetric,    $M$ contains three arcs $(2n+1-j_3, 2n-k)$, $(2n+1-j_2, 2n+1-k)$ and $(2n+1-j_1, 2n+2-k)$. According to the construction of $\phi$, $T_{2n+1-k }$ is obtained from $T_{2n-k+2}$ by inserting   $2n+1-j_1$, $T_{2n-k}$ is obtained from $T_{2n-k+1}$ by inserting   $2n+1-j_2$,  and $T_{2n-k-1}$ is obtained from $T_{2n-k}$ by inserting   $2n+1-j_3$.   Consider the insertion path $P(T_{2n+1-k}\leftarrow 2n+1-j_2)$ and the insertion path $P((T_{2n+1-k}\leftarrow 2n+1-j_2)\leftarrow 2n+1-j_3)$. As $j_2<j_3$, the insertion path $P((T_{2n+1-k}\leftarrow 2n+1-j_2)\leftarrow 2n+1-j_3)$ lies weakly  to the left of the insertion path $P(T_{2n+1-k}\leftarrow 2n+1-j_2)$.  Notice that each $T_i$ has at most two rows.    we have   (\rmnum{1}) $\lambda^{2n-k}$ is obtained from $\lambda^{2n+1-k}$ by adding  a square at the first row,  and (\rmnum{2})   $\lambda^{2n-k-1}$ is obtained from $\lambda^{2n-k}$ by adding a square at the second row.
    As $\lambda^i=\lambda^{2n-i}$,   we deduce that $k\in \mathrm{Val}(\phi(M))$, implying that $ \mathrm{Val}(M)\subseteq \mathrm{Val}(\phi(M))$.

    Now we proceed to show that $\mathrm{Val}(\phi(M))\subseteq \mathrm{Val}(M)$.
    Assume that $k\in \mathrm{Val}(\phi(M))$. As $\lambda^i=\lambda^{2n-i}$, we have (\rmnum{1}) $\lambda^{2n+1-k}$ is obtained from $\lambda^{2n+2-k}$ by adding  a square, (\rmnum{2}) $\lambda^{2n-k}$ is obtained from $\lambda^{2n+1-k}$ by adding  a square at the first row,  and (\rmnum{3})   $\lambda^{2n-k-1}$ is obtained from $\lambda^{2n-k}$ by adding a square at the second row. Assume that $T_{2n+1-k }$ is obtained from $T_{2n-k+2}$ by inserting   $j_1$, $T_{2n-k}$ is obtained from $T_{2n-k+1}$ by inserting   $j_2$,  and $T_{2n-k-1}$ is obtained from $T_{2n-k}$ by inserting   $j_3$. Consider the insertion path $P(T_{2n+2-k}\leftarrow j_1)$ and the insertion path $P((T_{2n+2-k}\leftarrow j_1)\leftarrow j_2)$.  We must have $j_1<j_2$. If not, then  the insertion path $P((T_{2n+2-k}\leftarrow j_1)\leftarrow j_2)$  lies weakly to the left of the insertion path $P(T_{2n+2-k}\leftarrow j_1)$.
    Then the insertion of $j_2$ in $T_{2n+1-k}$ would increase a square at the second row of $T_{2n+1-k}$.
     This contradicts   the fact   $\lambda^{2n-k}$ is obtained from $\lambda^{2n+1-k}$ by adding  a square at the first row. As $j_2$  is positioned at the end of the first row of $T_{2n-k}$,  if $j_3>j_2$, then when insert $j_3$ to $T_{2n-k}$, we just put $j_3$ at the end of the first row. This yields a contradiction with the fact that $\lambda^{2n-k-1}$ is obtained from $\lambda^{2n-k}$ by adding a square at the second row. Hence, we have $j_2>j_3$. Then $M$ contains three arcs $(k-1, 2n+1-j_1)$, $(k, 2n+1-j_2)$ and $(k+1, 2n+1-j_3)$ which form a valley of $M$.  Hence, we have $k\in \mathrm{Val}(M)$, implying that $\mathrm{Val}(\phi(M))\subseteq \mathrm{Val}(M)$.  Consequently,   we have concluded that $ \mathrm{Val}(M)= \mathrm{Val}(\phi(M))$ as desired, completing the proof. \qed

    By the same reasoning as in the proof of Lemma \ref{lemphi1}, one can deduce the following result.

\begin{lemma}\label{lemphi2}
The bijection $\phi$ induces a bijection between  $\mathcal{SNM}_3(n)$ and $\mathcal{SOC}_2(n)$ such that for a matching  $M\in \mathcal{SNM}_3(n)$, we have $\mathrm{type}(M)=\mathrm{type}(\phi(M))$ and $ \mathrm{Val}(M)= \mathrm{Peak}(\phi(M))$.

\end{lemma}

Recall that a {\em Dyck path} of length  $2n$ is a
	lattice path from  the origin $(0,0)$ to the destination $(2n, 0)$,   using only
	up steps $(1, 1)$ and down steps $(1, -1)$ and never  passing below the $x$-axis.    A  Dyck path of length $2n$  is said to be {\em symmetric  } if   its reflection about the
line  $x=n$ is itself.
A pair $(P,Q)$ of lattice  paths  is said  to be { \em noncrossing }  if $P$ never goes below $Q$.  Let $\mathcal{NC}(n)$  and  $\mathcal{SNC}_{n}$  denote the sets of  pairs of   noncrossing     Dyck paths of length $2n$ and  pairs of   noncrossing  symmetric   Dyck paths of length $2n$, respectively.

  Both oscillating  tableaux in $\mathcal{OR}_2(n)$  and oscillating  tableaux in $\mathcal{OC}_2(n)$   are proved  to be in one-to-one correspondence with pairs of noncrossing   Dyck paths in $\mathcal{NC}(n)$ by Chen-Deng-Du-Stanley-Yan \cite{ChenD}.

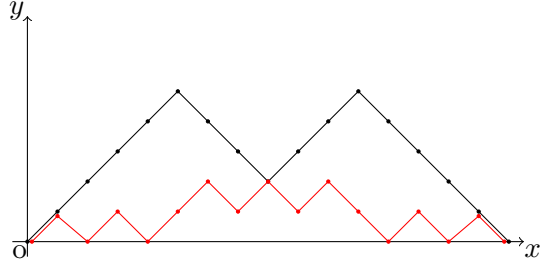
\begin{figure} [H]
  \centering
\begin{tikzpicture}[font =\small , scale = 0.4]
	\draw [->](-8,0)--(9,0);
	\draw [->](-7.5,-0.5)--(-7.5,7.5);
	\node at(-7.75,-0.35) {\small o};
	\node at(-7.85,7.7) {\small $y$};
	\node at(9.3,-0.3) {\small $x$};
	\filldraw[black](-7.5,0)circle(1.5pt);
	\draw (-7.5,0)--(-6.5,1);
	\filldraw[black](-6.5,1)circle(1.5pt);
	\draw (-6.5,1)--(-5.5,2);
	\filldraw[black](-5.5,2)circle(1.5pt);
	\draw (-5.5,2)--(-4.5,3);
	\filldraw[black](-4.5,3)circle(1.5pt);
	\draw (-4.5,3)--(-3.5,4);
	\filldraw[black](-3.5,4)circle(1.5pt);
	\draw (-3.5,4)--(-2.5,5);
	\filldraw[black](-2.5,5)circle(1.5pt);
	\draw (-2.5,5)--(-1.5,4);
	\filldraw[black](-1.5,4)circle(1.5pt);
	\draw (-1.5,4)--(-0.5,3);
	\filldraw[black](-0.5,3)circle(1.5pt);
	\draw (-0.5,3)--(0.5,2);
	\filldraw[black](0.5,2)circle(1.5pt);
	\draw (0.5,2)--(1.5,3);
	\filldraw[black](1.5,3)circle(1.5pt);
	\draw (1.5,3)--(2.5,4);
	\filldraw[black](2.5,4)circle(1.5pt);
	\draw (2.5,4)--(3.5,5);
	\filldraw[black](3.5,5)circle(1.5pt);
	\draw (3.5,5)--(4.5,4);
	\filldraw[black](4.5,4)circle(1.5pt);
	\draw (4.5,4)--(5.5,3);
	\filldraw[black](5.5,3)circle(1.5pt);
	\draw (5.5,3)--(6.5,2);
	\filldraw[black](6.5,2)circle(1.5pt);
	\draw (6.5,2)--(7.5,1);
	\filldraw[black](7.5,1)circle(1.5pt);
	\draw (7.5,1)--(8.5,0);
	\filldraw[black](8.5,0)circle(1.5pt);
	
	\filldraw[red](-7.35,0)circle(1.5pt);
	\draw[red] (-7.35,0)--(-6.5,0.85);
	\filldraw[red](-6.5,0.85)circle(1.5pt);
	\draw[red] (-6.5,0.85)--(-5.5,0);
	\filldraw[red](-5.5,0)circle(1.5pt);
	\draw[red] (-5.5,0)--(-4.5,1);
	\filldraw[red](-4.5,1)circle(1.5pt);
	\draw[red] (-4.5,1)--(-3.5,0);
	\filldraw[red](-3.5,0)circle(1.5pt);
	\draw[red] (-3.5,0)--(-2.5,1);
	\filldraw[red](-2.5,1)circle(1.5pt);
	\draw[red] (-2.5,1)--(-1.5,2);
	\filldraw[red](-1.5,2)circle(1.5pt);
	\draw[red] (-1.5,2)--(-0.5,1);
	\filldraw[red](-0.5,1)circle(1.5pt);
	\draw[red] (-0.5,1)--(0.5,2);
	\filldraw[red](0.5,2)circle(1.5pt);
	\draw[red] (0.5,2)--(1.5,1);
	\filldraw[red](1.5,1)circle(1.5pt);
	\draw[red] (1.5,1)--(2.5,2);
	\filldraw[red](2.5,2)circle(1.5pt);
	\draw[red] (2.5,2)--(3.5,1);
	\filldraw[red](3.5,1)circle(1.5pt);
	\draw[red] (3.5,1)--(4.5,0);
	\filldraw[red](4.5,0)circle(1.5pt);
	\draw[red] (4.5,0)--(5.5,1);
	\filldraw[red](5.5,1)circle(1.5pt);
	\draw[red] (5.5,1)--(6.5,0);
	\filldraw[red](6.5,0)circle(1.5pt);
	\draw[red] (6.5,0)--(7.5,0.85);
	\filldraw[red](7.5,0.85)circle(1.5pt);
	\draw[red] (7.5,0.85)--(8.35,0);
	\filldraw[red](8.35,0)circle(1.5pt);						
	\end{tikzpicture}
  \caption{An example of the bijection $\psi$.}\label{fig:OR2}
\end{figure}

\noindent {\bf The Bijection $\psi$ between $\mathcal{OR}_2(n)$ and $\mathcal{NC}(n)$}\\
Let $O=\{\lambda^0, \lambda^1, \ldots, \lambda^{2n}\}\in \mathcal{OR}_2(n)$ with $\lambda^{i}=(x_i, y_i)$.  Then,  a pair $\psi(O)=(P,Q)$ of noncrossing Dyck  paths is obtained  by letting $P=\{(i, x_i+y_i)\mid i=0, 1, \ldots, 2n\}$ and $Q=\{(i, x_i-y_i)\mid i=0, 1, \ldots, 2n\}$.
For example, let $O = (\emptyset,(1), (1, 1), (2, 1), (2, 2),$ $(3, 2), (3, 1), (2, 1), (2), (2, 1), (3, 1), (3, 2), (2, 2), (2, 1), (1, 1), (1), \emptyset)$ be an oscillating tableau in $\mathcal{OR}_2(8)$.
By using the map $\psi$, we obtain a pair of noncrossing Dyck paths as shown in Figure \ref{fig:OR2}.

\begin{figure}[H]
  \centering
\begin{center}
	\begin{tikzpicture}[font =\small , scale = 0.4]
		\draw [->](-8,0)--(9,0);
		\draw [->](-7.5,-0.5)--(-7.5,7.5);
		\node at(-7.75,-0.35) {\small o};
		\node at(-7.85,7.7) {\small $y$};
		\node at(9.3,-0.3) {\small $x$};
		\filldraw[black](-7.5,0)circle(1.5pt);
		\draw (-7.5,0)--(-6.5,1);
		\filldraw[black](-6.5,1)circle(1.5pt);
		\draw (-6.5,1)--(-5.5,2);
		\filldraw[black](-5.5,2)circle(1.5pt);
		\draw (-5.5,2)--(-4.5,3);
		\filldraw[black](-4.5,3)circle(1.5pt);
		\draw (-4.5,3)--(-3.5,4);
		\filldraw[black](-3.5,4)circle(1.5pt);
		\draw (-3.5,4)--(-2.5,5);
		\filldraw[black](-2.5,5)circle(1.5pt);
		\draw (-2.5,5)--(-1.5,4);
		\filldraw[black](-1.5,4)circle(1.5pt);
		\draw (-1.5,4)--(-0.5,3);
		\filldraw[black](-0.5,3)circle(1.5pt);
		\draw (-0.5,3)--(0.5,2);
		\filldraw[black](0.5,2)circle(1.5pt);
		\draw (0.5,2)--(1.5,3);
		\filldraw[black](1.5,3)circle(1.5pt);
		\draw (1.5,3)--(2.5,4);
		\filldraw[black](2.5,4)circle(1.5pt);
		\draw (2.5,4)--(3.5,5);
		\filldraw[black](3.5,5)circle(1.5pt);
		\draw (3.5,5)--(4.5,4);
		\filldraw[black](4.5,4)circle(1.5pt);
		\draw (4.5,4)--(5.5,3);
		\filldraw[black](5.5,3)circle(1.5pt);
		\draw (5.5,3)--(6.5,2);
		\filldraw[black](6.5,2)circle(1.5pt);
		\draw (6.5,2)--(7.5,1);
		\filldraw[black](7.5,1)circle(1.5pt);
		\draw (7.5,1)--(8.5,0);
		\filldraw[black](8.5,0)circle(1.5pt);
		
		\filldraw[red](-7.35,0)circle(1.5pt);
		\draw[red] (-7.35,0)--(-6.5,0.85);
		\filldraw[red](-6.5,0.85)circle(1.5pt);
		\draw[red] (-6.5,0.85)--(-5.5,1.85);
		\filldraw[red](-5.5,1.85)circle(1.5pt);
		\draw[red] (-5.5,1.85)--(-4.5,0.95);
		\filldraw[red](-4.5,0.95)circle(1.5pt);
		\draw[red] (-4.5,0.95)--(-3.5,2);
		\filldraw[red](-3.5,2)circle(1.5pt);
		\draw[red] (-3.5,2)--(-2.5,1);
		\filldraw[red](-2.5,1)circle(1.5pt);
		\draw[red] (-2.5,1)--(-1.5,2);
		\filldraw[red](-1.5,2)circle(1.5pt);
		\draw[red] (-1.5,2)--(-0.6,2.9);
		\filldraw[red](-0.6,2.9)circle(1.5pt);
		\draw[red] (-0.6,2.9)--(0.5,1.8);
		\filldraw[red](0.5,1.8)circle(1.5pt);
		\draw[red] (0.5,1.8)--(1.6,2.9);
		\filldraw[red](1.6,2.9)circle(1.5pt);
		\draw[red] (1.6,2.9)--(2.5,2);
		\filldraw[red](2.5,2)circle(1.5pt);
		\draw[red] (2.5,2)--(3.5,1);
		\filldraw[red](3.5,1)circle(1.5pt);
		\draw[red] (3.5,1)--(4.5,2);
		\filldraw[red](4.5,2)circle(1.5pt);
		\draw[red] (4.5,2)--(5.5,1);
		\filldraw[red](5.5,1)circle(1.5pt);
		\draw[red] (5.5,1)--(6.5,1.85);
		\filldraw[red](6.5,1.85)circle(1.5pt);
		\draw[red] (6.5,1.85)--(7.5,0.85);
		\filldraw[red](7.5,0.85)circle(1.5pt);
		\draw[red] (7.5,0.85)--(8.35,0);
		\filldraw[red](8.35,0)circle(1.5pt);
	\end{tikzpicture}
\end{center}

  \caption{An example of the bijection $\bar{\psi}$. }\label{fig:OC2}
\end{figure}
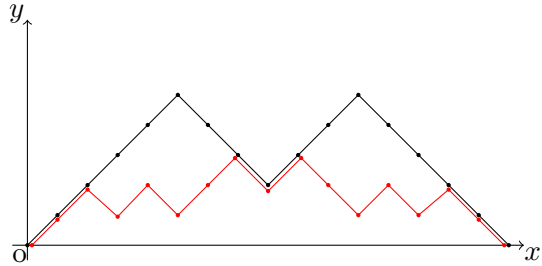

  \noindent {\bf The Bijection $\bar{\psi}$ between $\mathcal{OC}_2(n)$ and $\mathcal{NC}(n)$}\\
Let $O=\{\lambda^0, \lambda^1, \ldots, \lambda^{2n}\}\in \mathcal{OC}_2(n)$ with $(\lambda^{i})^{T}=(x_i, y_i)$.  Then,  a pair $\bar{\psi}(O)=(P,Q)$ of noncrossing Dyck  paths is obtained  by letting $P=\{(i, x_i+y_i)\mid i=0, 1, \ldots, 2n\}$ and $Q=\{(i, x_i-y_i)\mid i=0, 1, \ldots,  2n\}$.
For example, let $O = (\emptyset,(1), (1, 1), (2, 1), (2, 1, 1),$ $ (2, 2, 1), (2, 1, 1), (1, 1, 1), (1, 1), (1, 1, 1), (2, 1, 1), (2, 2, 1), (2, 1, 1), (2, 1), (1, 1), (1), \emptyset)$ be an oscillating tableau in $\mathcal{OC}_2(8)$.
By using the map $\bar{\psi}$, we obtain a pair of noncrossing Dyck paths as illustrated  in Figure \ref{fig:OC2}.

 For a lattice path $P$ of length $n$, denote by $P^{ w }$  the word obtained from $P$ by encoding each up step by the letter $U$ and each down step by  $D$.
Let $(P,Q)$ be a pair  of noncrossing    lattice paths of length $n$ with  $P^{w}=p_1p_2\cdots p_{n}$ and  $Q^{w}=q_1q_2\cdots q_{n}$.   Define
$$
\mathcal{A}(P,Q)=\{i\mid p_{i-1}p_{i}p_{i+1}=UUU, q_{i}q_{i+1}=UD,  1<i< n \}
$$
and
$$
\mathcal{B}(P,Q)=\{i\mid p_{i-1}p_{i}p_{i+1}=UUU, q_{i-1}q_{i}=UD,   1<i< n \}.
$$
For example, let $(P, Q)$ be the pair of noncrossing Dyck paths as shown in Figure \ref{fig:OC2}, we have $\mathcal{A}(P, Q)= \{2, 4\}$ and $\mathcal{B}(P, Q)=\{3, 10\}$.

From the construction of the bijections $\psi$ and $\bar{\psi}$, one can easily verify  the following results.

\begin{lemma}\label{lempsi1}
The bijection $\psi$ induces a bijection between  $\mathcal{SOR}_2(n)$ and $\mathcal{SNC}(n)$. Moreover,  for any  oscillating  tableau $O\in \mathcal{SOR}_2(n)$,  its corresponding  pair $\psi(O)=(P, Q)$ of noncrossing Dyck paths  verifies that $\mathrm{type}(O)=P^w$ and $ \mathrm{Val}(O)= \mathcal{A}(\psi(O))$.
\end{lemma}

\begin{lemma}\label{lempsi2}
The bijection $\bar{\psi}$ induces a bijection between  $\mathcal{SOC}_2(n)$ and $\mathcal{SNC}(n)$.  Moreover,  for any  oscillating  tableau $O\in \mathcal{SOC}_2(n)$,  its  corresponding  pair $\bar{\psi}(O)=(P, Q)$ of  noncrossing Dyck paths  verifies that $\mathrm{type}(O)=P^w$ and $ \mathrm{Peak}(O)= \mathcal{B}(\bar{\psi}(O))$.
\end{lemma}

 In the following, we aim to establish a bijection  $\theta:\mathcal{SNC}(n) \rightarrow \mathcal{SNC}(n)  $.
 \begin{lemma}\label{lemtheta}
  There exists a bijection $\theta:\mathcal{SNC}(n) \rightarrow \mathcal{SNC}(n)  $ such that for any pair $(P,Q)\in \mathcal{SNC}(n) $,  we have  $ \theta((P, Q))=(P,Q')$ with $\mathcal{A}(P,Q)=\mathcal{B}(P,Q')$.
 \end{lemma}

 Before we establish the bijection $\theta$, we define two transformations which will play essential roles in the construction of $\theta$.

\noindent{\bf The Transformation $\alpha$ } \\
Let $(S,R)$ be a pair  of noncrossing lattice   paths  where $S=\{(k, p_k)\mid i\leq k\leq j\}$   with $p_{k+1}-p_k=1$ for all $i\leq k < j$ and
 $R=\{(k, q_k)\mid i\leq k\leq j\}$ with $q_k\geq 0$ for all   $i\leq k \leq j$.   Let $R^{w}=w_1w_2\cdots w_{j-i}$.   We construct  a  lattice path $\alpha(R)=R'=\{(k, q'_k)\mid i\leq k\leq j\}$ by considering the following three cases.

 \noindent {\bf Case \upshape (\rmnum{1})} $w_1=D$. In this case,      let $R'$ be a lattice path  with the same origin as that of $R$ such that $(R')^w$ is obtained from $R^w$ by cyclicly shifting each letter one unit to its left, namely, $(R')^w= w_2\cdots w_{j-i}w_1$.    More precisely,  we have $q'_i=q_i,  q'_j=q_j$,  and $q'_k=q_{k+1}+1$ for all $i<k<j$.
See  Figure \ref{fig:alpha1} for an  illustration.

\begin{figure}[H]
  \centering
	\begin{tikzpicture}[font =\small , scale = 0.4]
		\draw[->] (-8,0)--(-2,0);
		\draw[->] (-8,0)--(-8,6);
		\node at(-8.25,-0.25) {\small o};
		\node at(-1.7,-0.25) {\small $x$};
		\node at(-8.25,6.2) {\small $y$};
		\draw[->] (-1,2.5)--(0.5,2.5);
		\node at(-0.3,3){\small $\alpha$};
		\draw[->] (2,0)--(8,0);
		\draw[->] (2,0)--(2,6);
		\node at(1.75,-0.25) {\small o};
		\node at(8.3,-0.25) {\small $x$};
		\node at(1.75,6.2) {\small $y$};
		
		\draw[dashed](-7,0)--(-7,5);
		\draw[dashed](-6,0)--(-6,5);
		\draw[dashed](-5,0)--(-5,5);
		\draw[dashed](-4,0)--(-4,5);
		\draw[dashed](-3,0)--(-3,5);
		\draw[dashed](-8,1)--(-3,1);
		\draw[dashed](-8,2)--(-3,2);
		\draw[dashed](-8,3)--(-3,3);
		\draw[dashed](-8,4)--(-3,4);
		\draw[dashed](-8,5)--(-3,5);
		
		\draw[dashed](3,0)--(3,5);
		\draw[dashed](4,0)--(4,5);
		\draw[dashed](5,0)--(5,5);
		\draw[dashed](6,0)--(6,5);
		\draw[dashed](7,0)--(7,5);
		\draw[dashed](2,1)--(7,1);
		\draw[dashed](2,2)--(7,2);
		\draw[dashed](2,3)--(7,3);
		\draw[dashed](2,4)--(7,4);
		\draw[dashed](2,5)--(7,5);
		
		\draw(-7,2)--(-4,5);
		\filldraw[black](-7,2)circle(1.5pt);
		\filldraw[black](-6,3)circle(1.5pt);
		\filldraw[black](-5,4)circle(1.5pt);
		\filldraw[black](-4,5)circle(1.5pt);
		\draw[red](-6.9,1.9)--(-6,1);
		\draw[red](-6,1)--(-5,2);
		\draw[red](-5,2)--(-4,1);
		\filldraw[red](-6.9,1.9)circle(1.5pt);
		\filldraw[red](-6,1)circle(1.5pt);
		\filldraw[red](-5,2)circle(1.5pt);
		\filldraw[red](-4,1)circle(1.5pt);
		
		\draw(3,2)--(6,5);
		\filldraw[black](3,2)circle(1.5pt);
		\filldraw[black](4,3)circle(1.5pt);
		\filldraw[black](5,4)circle(1.5pt);
		\filldraw[black](6,5)circle(1.5pt);
        \draw[red](3,1.87)--(4,2.87);
		\draw[red](4,2.87)--(5,2);
		\draw[red](5,2)--(6,1);
		\filldraw[red](3,1.87)circle(1.5pt);
		\filldraw[red](4,2.87)circle(1.5pt);
		\filldraw[red](5,2)circle(1.5pt);
		\filldraw[red](6,1)circle(1.5pt);
		
		%\node at(-5.2,-1) {  $(a)$};
		%\node at(4.7,-1) {\large $(b)$};		
	\end{tikzpicture}
  \caption{An example for  Case \upshape (\rmnum{1}) of $\alpha$.}\label{fig:alpha1}
\end{figure}

     \noindent {\bf Case \upshape (\rmnum{2})}   $w_1=U$ and $q_k>0$ for all $i<k\leq j$.    In this case,  let $R'$ be a lattice path  with the same origin as that of $R$ such that $(R')^w$ is obtained from $R^w$ by cyclicly shifting each letter one unit to its left, namely, $(R')^w= w_2\cdots w_{j-i}w_1$.   More precisely,  we have $q'_i=q_i,  q'_j=q_j$,  and $q'_k=q_{k+1}-1$ for all $i<k<j$.
      See Figure \ref{fig:alpha2}  for an illustration.
\begin{figure}[H]
  \centering
\begin{tikzpicture} [font = \small, scale = 0.4]
		
		\draw[->] (-8,0)--(-2,0);
		\draw[->] (-8,0)--(-8,6);
		\node at(-8.25,-0.25) {\small o};
		\node at(-1.7,-0.25) {\small $x$};
		\node at(-8.25,6.2) {\small$y$};
		\draw[->] (-1,2.5)--(0.5,2.5);
		\node at(-0.3,3){\small $\alpha$};
		\draw[->] (2,0)--(8,0);
		\draw[->] (2,0)--(2,6);
		\node at(1.75,-0.25) {\small o};
		\node at(8.3,-0.25) {\small $x$};
		\node at(1.75,6.2) {\small $y$};
		
		\draw[dashed](-7,0)--(-7,5);
		\draw[dashed](-6,0)--(-6,5);
		\draw[dashed](-5,0)--(-5,5);
		\draw[dashed](-4,0)--(-4,5);
		\draw[dashed](-3,0)--(-3,5);
		\draw[dashed](-8,1)--(-3,1);
		\draw[dashed](-8,2)--(-3,2);
		\draw[dashed](-8,3)--(-3,3);
		\draw[dashed](-8,4)--(-3,4);
		\draw[dashed](-8,5)--(-3,5);
		
		\draw[dashed](3,0)--(3,5);
		\draw[dashed](4,0)--(4,5);
		\draw[dashed](5,0)--(5,5);
		\draw[dashed](6,0)--(6,5);
		\draw[dashed](7,0)--(7,5);
		\draw[dashed](2,1)--(7,1);
		\draw[dashed](2,2)--(7,2);
		\draw[dashed](2,3)--(7,3);
		\draw[dashed](2,4)--(7,4);
		\draw[dashed](2,5)--(7,5);
		
		\draw(-7,2)--(-4,5);
		\filldraw[black](-7,2)circle(1.5pt);
		\filldraw[black](-6,3)circle(1.5pt);
		\filldraw[black](-5,4)circle(1.5pt);
		\filldraw[black](-4,5)circle(1.5pt);
		\draw[red](-7,1)--(-5,3);
		\draw[red](-5,3)--(-4,2);
		\filldraw[red](-7,1)circle(1.5pt);
		\filldraw[red](-6,2)circle(1.5pt);
		\filldraw[red](-5,3)circle(1.5pt);
		\filldraw[red](-4,2)circle(1.5pt);
		
		\draw(3,2)--(6,5);
		\filldraw[black](3,2)circle(1.5pt);
		\filldraw[black](4,3)circle(1.5pt);
		\filldraw[black](5,4)circle(1.5pt);
		\filldraw[black](6,5)circle(1.5pt);
		\draw[red](3,1)--(4,2);
		\draw[red](4,2)--(5,1);
		\draw[red](5,1)--(6,2);
		\filldraw[red](3,1)circle(1.5pt);
		\filldraw[red](4,2)circle(1.5pt);
		\filldraw[red](5,1)circle(1.5pt);
		\filldraw[red](6,2)circle(1.5pt);
		
		%\node at(-5.2,-1) {\large $(a)$};
		%\node at(4.7,-1) {\large $(b)$};		
	\end{tikzpicture}
  \caption{An example for  Case \upshape (\rmnum{2}) of $\alpha$.}\label{fig:alpha2}
\end{figure}

                  \noindent {\bf Case \upshape (\rmnum{3})}  $w_1=U$ and there exists some integer $i<k\leq j$ such  that $q_k=0$. Let $\ell$ be  the smallest such  integer.  Clearly,  $w_{\ell-i}=D$ according to the selection of $\ell$. Then, let
                   $R'$  be a lattice path with the same  origin as that of $R$  such that $(R')^w$ is obtained from $R^w$ by replacing $w_1$ by $D$,  replacing $w_{\ell-i}$ by $U$, and then cyclicly shifting each letter one unit to its left, namely,
                    $(R')^w=w_2\cdots w_{\ell-i-1}U w_{\ell-i+1}\cdots w_{j-i} D.$
                  More precisely,  we have $q'_i=q_i,  q'_j=q_j$, $q'_k=q_{k+1}-1$ for all $i<k<\ell-1$,  and $q'_k=q_{k+1}+1$  for all $\ell-1\leq k<j$.
                   See Figure \ref{fig:alpha3} for an illustration.

\begin{figure}[H]
  \centering
\begin{tikzpicture}[font = \small, scale = 0.4]	
		\draw[->] (-8,0)--(-2,0);
		\draw[->] (-8,0)--(-8,6);
		\node at(-8.25,-0.25) {\small o};
		\node at(-1.7,-0.25) {\small $x$};
		\node at(-8.25,6.2) {\small $y$};
		\draw[->] (-1,2.5)--(0.5,2.5);
		\node at(-0.3,3){\small $\alpha$};
		\draw[->] (2,0)--(8,0);
		\draw[->] (2,0)--(2,6);
		\node at(1.75,-0.25) {\small o};
		\node at(8.3,-0.25) {\small $x$};
		\node at(1.75,6.2) {\small $y$};
		
		\draw[dashed](-7,0)--(-7,5);
		\draw[dashed](-6,0)--(-6,5);
		\draw[dashed](-5,0)--(-5,5);
		\draw[dashed](-4,0)--(-4,5);
		\draw[dashed](-3,0)--(-3,5);
		\draw[dashed](-8,1)--(-3,1);
		\draw[dashed](-8,2)--(-3,2);
		\draw[dashed](-8,3)--(-3,3);
		\draw[dashed](-8,4)--(-3,4);
		\draw[dashed](-8,5)--(-3,5);
		
		\draw[dashed](3,0)--(3,5);
		\draw[dashed](4,0)--(4,5);
		\draw[dashed](5,0)--(5,5);
		\draw[dashed](6,0)--(6,5);
		\draw[dashed](7,0)--(7,5);
		\draw[dashed](2,1)--(7,1);
		\draw[dashed](2,2)--(7,2);
		\draw[dashed](2,3)--(7,3);
		\draw[dashed](2,4)--(7,4);
		\draw[dashed](2,5)--(7,5);
		
		\draw(-8,0)--(-3,5);
		\filldraw[black](-8,0)circle(1.5pt);
		\filldraw[black](-7,1)circle(1.5pt);
		\filldraw[black](-6,2)circle(1.5pt);
		\filldraw[black](-5,3)circle(1.5pt);
		\filldraw[black](-4,4)circle(1.5pt);
		\filldraw[black](-3,5)circle(1.5pt);
		\draw[red](-7.87,0)--(-7,0.87);
		\draw[red](-7,0.87)--(-6,0);
		\draw[red](-6,0)--(-5,1);
		\draw[red](-5,1)--(-4,0);
		\draw[red](-4,0)--(-3,1);
		\filldraw[red](-7.87,0)circle(1.5pt);
		\filldraw[red](-7,0.87)circle(1.5pt);
		\filldraw[red](-6,0)circle(1.5pt);
		\filldraw[red](-5,1)circle(1.5pt);
		\filldraw[red](-4,0)circle(1.5pt);
		\filldraw[red](-3,1)circle(1.5pt);
		
		\draw(2,0)--(7,5);
		\filldraw[black](2,0)circle(1.5pt);
		\filldraw[black](3,1)circle(1.5pt);
		\filldraw[black](4,2)circle(1.5pt);
		\filldraw[black](5,3)circle(1.5pt);
		\filldraw[black](6,4)circle(1.5pt);
		\filldraw[black](7,5)circle(1.5pt);
		\draw[red](2.13,0)--(3,0.87);
		\draw[red](3,0.87)--(4,1.87);
		\draw[red](4,1.87)--(5,1);
		\draw[red](5,1)--(6,2);
		\draw[red](6,2)--(7,1);
		\filldraw[red](2.13,0)circle(1.5pt);
		\filldraw[red](3,0.87)circle(1.5pt);
		\filldraw[red](4,1.87)circle(1.5pt);
		\filldraw[red](5,1)circle(1.5pt);
		\filldraw[red](6,2)circle(1.5pt);
		\filldraw[red](7,1)circle(1.5pt);
		
		%\node at(-5.2,-1) {\large $(a)$};
		%\node at(4.7,-1) {\large $(b)$};				
	\end{tikzpicture}
  \caption{An example for  Case \upshape (\rmnum{3}) of $\alpha$.}\label{fig:alpha3}
\end{figure}

In the following, we aim to show that the transformation $\alpha$ has the following property.
 \begin{lemma}\label{lemalpha}
 The resulting pair $(S, R')$ is noncrossing   such that
   \begin{itemize}
     \item[(1)]  the origin (resp. the  destination) of $ R'$ is the same as that of $R$;
     \item[(2)] the lattice path $R'$  lies weakly above the $x$-axis;
     \item[(3)]  $\mathcal{A}(S, R )=\mathcal{B}(S, R')$.
      \end{itemize}
   \end{lemma}
\pf  Properties (1) and (2) follows directly from the construction of $\alpha$.
   In order to show that $(S, R')$ is noncrossing, it suffices  to show that $q'_k\leq p_k$ for all $i\leq k\leq j$.
 It is apparent that the assertion holds for $k=i$ or $k=j$ as $q'_i=q_i\leq p_i$ and $q'_j=q_j\leq p_j$.  Now we assume that $i<k<j$. we  have three cases.

 \noindent{\bf Case 1}  $w_1=D$. In this case,    we have   $q'_k=q_{k+1}+1$.  As $q_{i+1}-q_i=-1$ and $p_{t+1}-p_t=1$ for all $i\leq t< j$, we conclude that $p_{t}-q_t\geq 2$ for all $i< t\leq j$. This yields that $q'_k=q_{k+1}+1\leq p_{k+1}-1=p_k$ as desired.

\noindent {\bf Case 2}  $ w_1=U$ and $q_t>0$ for all $i<t\leq j$. In this case,   we have  $q'_k=q_{k+1}-1$.  Then  $p_{k}=p_{k+1}-1$ together with $q_{k+1}\leq p_{k+1}$ implies  that
$q'_k=q_{k+1}-1\leq p_{k+1}-1=p_k$ as desired.

                  \noindent {\bf Case 3}  $ w_1=U$ and there exists some integer $i<t\leq j$ such  that $q_t=0$. Let $\ell$ be  the smallest such  integer. Then we have  $q'_k=q_{k+1}-1$  when $i<k<\ell-1$  and $q'_k=q_{k+1}+1$ otherwise.  In the former case, $p_{k}=p_{k+1}-1$ together with $q_{k+1}\leq p_{k+1}$ implies  that $q'_k=q_{k+1}-1\leq p_{k+1}-1=p_k$ as desired.  In the latter case, we have $q_{k+1}\leq p_{k+1}-2$ since   $w_{\ell-i}=D$. Hence we have $q'_k=q_{k+1}+1\leq p_{k+1}-1=p_{k}$ as desired.

So far,  we have concluded that $(S, R')$ is noncrossing. Next we aim to show that $\mathcal{A}(S, R )=\mathcal{B}(S, R')$.   Let $(R')^w=v_1v_2\cdots v_{j-i}$.  We have two cases.

 \noindent{\bf Case 1}   $ w_1=U$ and $q_k>0$ for all $i<k\leq j$ or $w_1=D$. In this case,  we have $(R')^w= w_2\ldots w_{j-i}w_1$. One can easily check that
 $w_kw_{k+1}=UD$ if and only if $v_{k-1}v_k=UD$  for all $1<k<j-i$. This implies that  $\mathcal{A}(S, R )=\mathcal{B}(S, R')$.

  \noindent {\bf Case 2}  $ w_1=U$ and there exists some integer $i<t\leq j$ such  that $q_t=0$.  Let $\ell$ be  the smallest such  integer.
   In this case, we have
 $(R')^w=w_2\cdots w_{\ell-i-1}U w_{\ell-i+1}\cdots w_{j-i} D.$
   According to the selection of $\ell$, we must have either $w_{\ell-i-1}=D$ or $\ell=i+2$.  This implies that $\ell-i-1\notin \mathcal{A}(S, R)$. As $w_{\ell-i}=D$, we have $\ell-i\notin   \mathcal{A}(S, R)$.
     As $q_k\geq 0$ for all $i\leq k\leq j$ and $q_\ell=0$, we must have $w_{\ell-i+1}=U$ or $\ell=j$.  As    $v_{\ell-i-1}=U$ and $v_{\ell-i}=w_{\ell-i+1}=U$ when $\ell<j$, we have $\ell - i -1,\ell-i\notin \mathcal{B}(S, R') $.
     %Moreover, $v_{\ell-i-1}=U$ implies that $\ell-i-1\notin \mathcal{B}(S, R')$.
   It is not difficult to check that
 $w_kw_{k+1}=UD$ if and only if $v_{k-1}v_k=UD$  for all $1<k< j-i$ and $k\notin \{\ell-i-1, \ell-i\}$. This yields that $\mathcal{A}(S, R )=\mathcal{B}(S, R')$, completing the proof.
\qed

 \noindent{\bf The Transformation $\beta$.} \\
 Let $(S,R')$ be a pair  of noncrossing lattice   paths  where $S=\{(k, p_k)\mid i\leq k\leq j\}$   with $p_{k+1}-p_k=1$ for all $i\leq k<j$ and
 $R'=\{(k, q'_k)\mid i\leq k\leq j\}$ with $q'_k\geq 0$ for all all $i\leq k \leq j$.   Let $(R')^{w}=v_1v_2\cdots v_{j-i}$.  We construct  a  lattice path  $\beta(R')=R=\{(k, q_k)\mid i\leq k\leq j\}$ by considering the following three cases.

 \noindent {\bf Case \upshape (\rmnum{1}$'$)} $v_{j-i}=D$ and $q'_k>0$ for all $i\leq k<j$. In this case, let $R$ be a lattice path with the same origin as that of $R'$ such that $R^w$ is obtained from $(R')^w$ by cyclicly shifting each letter one unit to its right, namely,  $ R^w= v_{j-i}v_1v_2\cdots v_{j-i-1}$.    More precisely,  we have $q_i=q'_i,  q_j=q'_j$,  and $q_{k}=q'_{k-1}-1$ for all $i<k<j$.

     \noindent {\bf Case \upshape (\rmnum{2}$'$)}  $v_{j-i}=U$.   In this case, let $R$ be a lattice path with the same origin as that of $R'$ such that $R^w$ is obtained from $(R')^w$ by cyclicly shifting each letter one unit to its right, namely,  $(R)^w= v_{j-i}v_1v_2\cdots v_{j-i-1}$.  More precisely,  $q_i=q'_i,  q_j=q'_j$,  and $q_{k}=q'_{k-1}+1$ for all $i<k<j$.

                  \noindent {\bf Case \upshape (\rmnum{3}$'$)}  $v_{j-i}=D$ and there exists some integer $i\leq k< j$ such  that $q_k=0$. Let $\ell$ be  the greatest  such  integer.  Let $R$ be a lattice path with the same origin as that of $R'$  such that  $R^w$ is obtained from $(R')^w$ by  replacing $v_{j-i}$ by $U$, replacing $v_{\ell-i+1}$ by $D$, and then cyclicly shifting each letter one unit to its right,  namely,  $(R)^w= Uv_1v_2\cdots v_{\ell-i} D v_{\ell-i+2}\cdots  v_{j-i-1}$.
                  More precisely,  we have $q_i=q'_i,  q_j=q'_j$,  $q_{k+1}=q'_{k}+1$ for all $i \leq k \leq \ell$,  and $q_{k+1}=q'_{k}-1$  for all $\ell < k<j-1$.

The following property of $\beta$ can be deduced by the same reasoning as in the proof of Lemma \ref{lemalpha}, and the proof is omitted here.
\begin{lemma}\label{lembeta}
The resulting pair $(S, R)$ is noncrossing   such that
   \begin{itemize}
     \item[(1)]  the origin (resp. the  destination) of $R=\beta(R')$ is the same as that of $R'$;
     \item[(2)] the lattice path $R$ lying weakly above the $x$-axis.
          \end{itemize}
\end{lemma}

Let $P=\{(i, p_i)\mid i=0, 1, \ldots, n\}$   be a lattice path with up steps and down steps.
A subpath $S=\{(k, p_k)\mid i\leq k\leq j\}$ for some $0\leq i< j\leq n$ is said to be a {\em  maximal  chain} of  up steps of $P$ if
$p_k-p_{k-1}=1$ for all $i<k\leq j$, $p_{i}-p_{i-1}=-1$ when $i>0$, and $p_{j+1}-p_j=-1$ when $j<n$.
  Similarly, a subpath $S=\{(k, p_k)\mid i\leq k\leq j\}$ for some $0\leq i< j\leq n$ is said to be a {\em  maximal  chain} of  down steps of $P$ if
$p_k-p_{k-1}=-1$ for all $i<k\leq j$, $p_{i}-p_{i-1}=1$ when $i>0$, and $p_{j+1}-p_j=1$ when $j<n$.

Now we are ready for the proof of Lemma \ref{lemtheta}.

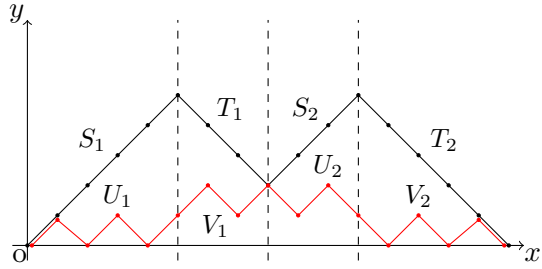
\begin{figure} [H]
\centering
	\begin{tikzpicture}[font =\small , scale = 0.4]
		\draw [->](-8,0)--(9,0);
		\draw [->](-7.5,-0.5)--(-7.5,7.5);
		\draw [dashed] (-2.5,-0.5)--(-2.5,7.5);
		\draw [dashed] (0.5,-0.5)--(0.5,7.5);
		\draw [dashed] (3.5,-0.5)--(3.5,7.5);
		\node at(-7.75,-0.35) {\small o};
		\node at(-7.85,7.7) {\small $y$};
		\node at(9.3,-0.3) {\small $x$};
		\node at(-5.35,3.5) {\footnotesize $S_1$};
		\node at(-4.5,1.65) {\footnotesize $U_1$};
		\node at(-0.75,4.55) {\footnotesize $T_1$};
		\node at(-1.25,0.65) {\footnotesize $V_1$};
		\node at(1.75,4.55) {\footnotesize $S_2$};
		\node at(2.5,2.65) {\footnotesize $U_2$};
		\node at(6.35,3.5) {\footnotesize $T_2$};
		\node at(5.5,1.65) {\footnotesize $V_2$};
		\filldraw[black](-7.5,0)circle(1.5pt);
		\draw (-7.5,0)--(-6.5,1);
		\filldraw[black](-6.5,1)circle(1.5pt);
		\draw (-6.5,1)--(-5.5,2);
		\filldraw[black](-5.5,2)circle(1.5pt);
		\draw (-5.5,2)--(-4.5,3);
		\filldraw[black](-4.5,3)circle(1.5pt);
		\draw (-4.5,3)--(-3.5,4);
		\filldraw[black](-3.5,4)circle(1.5pt);
		\draw (-3.5,4)--(-2.5,5);
		\filldraw[black](-2.5,5)circle(1.5pt);
		\draw (-2.5,5)--(-1.5,4);
		\filldraw[black](-1.5,4)circle(1.5pt);
		\draw (-1.5,4)--(-0.5,3);
		\filldraw[black](-0.5,3)circle(1.5pt);
		\draw (-0.5,3)--(0.5,2);
		\filldraw[black](0.5,2)circle(1.5pt);
		\draw (0.5,2)--(1.5,3);
		\filldraw[black](1.5,3)circle(1.5pt);
		\draw (1.5,3)--(2.5,4);
		\filldraw[black](2.5,4)circle(1.5pt);
		\draw (2.5,4)--(3.5,5);
		\filldraw[black](3.5,5)circle(1.5pt);
		\draw (3.5,5)--(4.5,4);
		\filldraw[black](4.5,4)circle(1.5pt);
		\draw (4.5,4)--(5.5,3);
		\filldraw[black](5.5,3)circle(1.5pt);
		\draw (5.5,3)--(6.5,2);
		\filldraw[black](6.5,2)circle(1.5pt);
		\draw (6.5,2)--(7.5,1);
		\filldraw[black](7.5,1)circle(1.5pt);
		\draw (7.5,1)--(8.5,0);
		\filldraw[black](8.5,0)circle(1.5pt);
		
		\filldraw[red](-7.35,0)circle(1.5pt);
		\draw[red] (-7.35,0)--(-6.5,0.85);
		\filldraw[red](-6.5,0.85)circle(1.5pt);
		\draw[red] (-6.5,0.85)--(-5.5,0);
		\filldraw[red](-5.5,0)circle(1.5pt);
		\draw[red] (-5.5,0)--(-4.5,1);
		\filldraw[red](-4.5,1)circle(1.5pt);
		\draw[red] (-4.5,1)--(-3.5,0);
		\filldraw[red](-3.5,0)circle(1.5pt);
		\draw[red] (-3.5,0)--(-2.5,1);
		\filldraw[red](-2.5,1)circle(1.5pt);
		\draw[red] (-2.5,1)--(-1.5,2);
		\filldraw[red](-1.5,2)circle(1.5pt);
		\draw[red] (-1.5,2)--(-0.5,1);
		\filldraw[red](-0.5,1)circle(1.5pt);
		\draw[red] (-0.5,1)--(0.5,2);
		\filldraw[red](0.5,2)circle(1.5pt);
		\draw[red] (0.5,2)--(1.5,1);
		\filldraw[red](1.5,1)circle(1.5pt);
		\draw[red] (1.5,1)--(2.5,2);
		\filldraw[red](2.5,2)circle(1.5pt);
		\draw[red] (2.5,2)--(3.5,1);
		\filldraw[red](3.5,1)circle(1.5pt);
		\draw[red] (3.5,1)--(4.5,0);
		\filldraw[red](4.5,0)circle(1.5pt);
		\draw[red] (4.5,0)--(5.5,1);
		\filldraw[red](5.5,1)circle(1.5pt);
		\draw[red] (5.5,1)--(6.5,0);
		\filldraw[red](6.5,0)circle(1.5pt);
		\draw[red] (6.5,0)--(7.5,0.85);
		\filldraw[red](7.5,0.85)circle(1.5pt);
		\draw[red] (7.5,0.85)--(8.35,0);
		\filldraw[red](8.35,0)circle(1.5pt);						
	\end{tikzpicture}
	\caption{The decomposition of the pair $(P, Q)$ of noncrossing Dyck paths.}\label{decom}
\end{figure}

\noindent  {\bf Proof of Lemma \ref{lemtheta}.}  First we give a description of the map $\theta: \mathcal{SNC}(n)\rightarrow \mathcal{SNC}(n)$. Let $(P, Q)\in \mathcal{SNC}(n)$. Then $P$ can be uniquely decomposed as
$
S_1T_1S_2$ $T_2\cdots S_kT_k
$
where each $S_i$ is a maximal chain of up steps of $P$ and each $T_i$ is a maximal chain of down steps of $P$. Accordingly, the path $Q$ can be uniquely decomposed as
$U_1V_1U_2V_2\cdots U_kV_k$, where the $x$-coordinate of the origin of  $U_i$ (resp. $V_i$) is the same as that of $S_i$ (resp. $T_i$) and the $x$-coordinate of the destination  of  $U_i$ (resp. $V_i$) is the same as that of $S_i$ (resp. $T_i$).  Then,  we construct a lattice path
$$
Q'=\alpha(U_1)V'_1\alpha(U_2)V'_2\cdots \alpha(U_k)V'_k,
$$
where $V'_i  $ is obtained from $\alpha(U_{k+1-i})$  by reflecting  about the line $x=n$. Define $\theta((P, Q))=(P, Q')$.  For example, let $(P, Q)$ be the pair of noncrossing Dyck paths of length $16$ as illustrated in Figure \ref{fig:OR2}.  Then $P$ and $Q$ can be uniquely decomposed as shown in Figure \ref{decom}.  By applying the map $\theta$ to $(P, Q)$, its corresponding pair $(P, Q')$ is shown in Figure \ref{fig:OC2}.

Lemma \ref{lemalpha} ensures that resulting pair $(P, Q')\in \mathcal{NC}(n)$. Moreover, by the construction of $V'_i$,  the path $Q'$ is   symmetric.  Hence, we have concluded that $\theta((P,Q))=(P, Q')\in \mathcal{SNC}(n)$.
 Notice that
$\mathcal{A}(P, Q)$ is uniquely determined by the sets $\mathcal{A}(S_i, U_i)$ for all $1\leq i\leq k$ and $\mathcal{B}(P, Q')$  is uniquely determined by the sets $\mathcal{B}(S_i, \alpha(U_i))$ for all $1\leq i\leq k$.   Again by Lemma \ref{lemalpha}, we have $\mathcal{A}(S_i, U_i)=\mathcal{B}(S_i, \alpha(U_i))$. Thus, we deduce that $\mathcal{A}(P, Q)=\mathcal{B}(P, Q')$ as desired.

In order to show the map $\theta$ is a bijection, we  establish a map $\theta^{'}: \mathcal{SNC}(n)\rightarrow \mathcal{SNC}(n)$. Let $(P, Q)\in \mathcal{SNC}(n)$. Then $P$ can be uniquely decomposed as
$
S_1T_1S_2T_2\cdots $ $S_kT_k$,
where each $S_i$ is a maximal chain of up steps of $P$ and each $T_i$ is a maximal chain of down steps of $P$. Accordingly, the path $Q$ can be uniquely decomposed as
$U_1V_1U_2V_2\cdots U_kV_k$, where the $x$-coordinate of the origin of  $U_i$ (resp. $V_i$) is the same as that of $S_i$ (resp. $T_i$) and the $x$-coordinate of the destination  of  $U_i$ (resp. $V_i$) is the same as that of $S_i$ (resp. $T_i$).  Then,  we construct a lattice path
$$
Q'=\beta(U_1)V'_1\beta(U_2)V'_2\cdots \beta(U_k)V'_k,
$$
where $V'_i  $ is obtained from $\beta(U_{k+1-i})$  by reflecting  about the line $x=n$. Set $\theta^{'}((P,Q))=(P, Q')$.  Lemma \ref{lembeta} together with the construction of $V'_i$ guarantees that the resulting pair $(P, Q')\in \mathcal{SNC}(n)$.

In order to show that the maps $\theta$ and $\theta^{'}$ are inverses of each other, it suffices to show that $\alpha$ and $\beta$ are inverses of each other.
This can be justified by the fact that  cases (\rmnum{1}$'$), (\rmnum{2}$'$) and (\rmnum{3}$'$) in the construction of  $\beta$
correspond, respectively, to cases (\rmnum{1}), (\rmnum{2}) and (\rmnum{3}) of the construction of  $\alpha$. This completes the proof. \qed

\noindent  {\bf Proof of Theorem \ref{mainth2}.}  Combining  Lemmas \ref{lemchi1},  \ref{lemchi2}, \ref{lemphi1}, \ref{lemphi2}, \ref{lempsi1}, \ref{lempsi2} and \ref{lemtheta}, the map $ \Psi=\chi'\circ \phi^{-1} \circ \bar{\psi}^{-1}\circ\theta \circ\psi \circ\phi \circ\chi$ serves a     bijection between $\mathcal{ST}_{\lambda}(J_3)$ and $\mathcal{ST}_{\lambda}(I_3)$ such that for any $T\in \mathcal{ST}_{\lambda}(J_3)$, we have $\mathrm{Peak}(T)=\mathrm{Peak}(\Psi(T))$ for any self-conjugate Young diagram $\lambda$.  This completes the proof. \qed

\subsection{Proof of Theorem  \ref{mainth1}}

In  this  section,  we  will  adapt  the   proof  of  Theorem    \ref{BWX1} in \cite{BWX}  to  prove Theorem \ref{mainth1}.

\noindent{\bf Proof of Theorem \ref{mainth1}.}   Let $\lambda$ be a self-conjugate Young  diagram $\lambda$ with $k$ columns.  We proceed to establish a bijection $\Phi: \mathcal{ST}_{\lambda }(J_3\oplus \tau)\rightarrow \mathcal{ST}_{\lambda}(I_3\oplus \tau)$. Given a transversal $T=\{(i,t_i)\}_{i=1}^{k}\in \mathcal{ST}_{\lambda}(J_3\oplus \tau)$,   colour the square $(c,r)$ of $\lambda$ by white if the board of $\lambda$ lying below and to the right of it contains either $\tau$ or $\tau^{-1}$, or gray otherwise.  Then find the  $1$'s coloured by gray, and colour the corresponding rows and columns gray. Denote the white board by $\lambda'$ and let $T'$ be the  $01$-filling of $ T$ located at the  the board $\lambda'$. Clearly, the board $\lambda'$ is a self-conjugate Young diagram and  $T'$ is a symmetric transversal of $ \lambda'$.    $T$ avoids the pattern $J_3\oplus \tau$  implies that $T'$ avoids the pattern $J_3$.   Theorem \ref{mainth2} ensures that   by applying the map $\Psi$ to $T'$, we can get a symmetric  transversal $L'=\Psi(T')$ of $\lambda'$ such that  $L'$ avoids the pattern $I_3$ and $\mathrm{Peak}(L')= \mathrm{Peak}(T')$. Restoring the gray cells of $\lambda$ and  their contents, we obtain a symmetric transversal $L$  of $\lambda$ which avoids the pattern $I_3\oplus \tau$.

 Now we proceed show that $\mathrm{Peak}(L)=\mathrm{Peak}(T)$.
It is easily seen that if $i\in \mathrm{Peak}(T)$ and column $i$ is coloured by gray, then we have $i\in \mathrm{Peak}(L)$.
Now we assume that $i\in \mathrm{Peak}(T)$ and the lowest coloured white square in column $i$ is $(i, r)$.  Assume that column $i$ of $\lambda$ is numbered by $i'$ in $\lambda'$.  As $t_{i-1}<t_i>t_{i+1}$ and $c_{i-1}(\lambda)=c_i(\lambda)=c_{i+1}(\lambda)$, the lowest coloured  white square in   column $i-1$  (resp. $i+1$)  is $(i-1, r)$ (resp. $(i+1, r)$). This implies that $c_{i'-1}(\lambda')=c_{i'}(\lambda')=c_{i'+1}(\lambda')$ and the relative positions of $1$'s in columns $i'-1$, $i'$ and $i'+1$ in $T'$ are the same as those of  $1$'s in columns $i-1$, $i$ and $i+1$ in $T$.  Hence we have $i'\in  \mathrm{Peak}(T') =\mathrm{Peak}(L')$ by Theorem \ref{mainth2}. According to the construction of $L$,   the relative positions of $1$'s in columns $i'-1$, $i'$ and $i'+1$ in $L'$ are the same as those of  $1$'s in columns $i-1$, $i$ and $i+1$ in $L$. This yields that
 $i\in \mathrm{Peak}(L)$. Hence, we have $\mathrm{Peak}(T) \subseteq \mathrm{Peak}(L)$.  By similar arguments, one can verify that  $\mathrm{Peak}(L) \subseteq \mathrm{Peak}(T)$. Therefore, we have concluded that $\mathrm{Peak}(L)=\mathrm{Peak}(T)$.

In order to show that the map   $\Phi$ is a bijection, it suffices to show that  the above procedure is invertible.  It is obvious that    $\Phi$ changes the $01$-filling located at the white squares  and leaves the $01$-filling located at the  gray  squares  fixed.   Hence  when  applying  the  inverse  map   $\Phi^{-1}$,   the  colouring  of $L$ will result in the same  Young diagram  $\lambda'$  and the same transversal $L'$  such that when  applying the inverse bijection   $\Psi^{-1}$ to $L'$, we will recover the same transversal $T'$ and hence the same transversal $T$.  This completes the proof. \qed

\noindent{\bf Proof of Theorem \ref{generalth}.} Note that the permutation matrix of an   involution  $\pi\in\mathcal{I}_n$ is a  symmetric  transversal of      $\lambda=(\lambda_1, \lambda_2, \ldots, \lambda_n)$ with $\lambda_1=\lambda_2=\cdots=\lambda_n=n$. Thus, by  Theorem \ref{mainth1},    the map $\Phi$  serves as   a bijection between    $\mathcal{I}_n(321 \oplus \tau)$ and $\mathcal{I}_{n}(123 \oplus \tau)$ such that  for any $\pi\in \mathcal{I}_n(321 \oplus \tau)$, we have $\mathrm{Peak}(\pi)=\mathrm{Peak}(\Phi(\pi))$. This completes the proof.

\noindent{\bf Proof of Conjecture \ref{conj3}}. Notice that the  permutation matrix of an alternating involution $\pi\in \mathcal{AI}_{n}$  is a symmetric transversal of a square Young diagram $\lambda=(\lambda_1, \lambda_2, \ldots, \lambda_n)$ with $\lambda_1=\lambda_2=\cdots=\lambda_2=n$  such that $\mathrm{Peak}(\pi)=\{2, 4, \ldots, 2\lfloor {n-1\over 2}\rfloor\}$.  Let $\pi\in \mathcal{AI}_n(321\oplus \tau) $  and $\Phi(\pi)=\sigma=\sigma_1\sigma_2\cdots \sigma_{n}$.
 In order to show that $\Phi(\pi)\in \mathcal{AI}_n(123\oplus \tau) $,      by Theorem \ref{generalth}, it suffices to show that  $\sigma_{2k}>\sigma_{2k-1}$  when $n=2k$.
 This is justified by the fact that, when applying the bijection $\Phi$,  column $2k$  is always coloured by gray and  all the squares coloured by white (if any) in column $2k-1$  are positioned  above the $1$ located in column $2k$ as  $\tau$ is nonempty.
 Hence $\Phi$ induces a bijection between $\mathcal{AI}_n(321\oplus \tau)$ and $\mathcal{AI}_n(123\oplus \tau)$.
 This completes the proof.
 \qed

	\section*{Acknowledgments}
	The work  was supported by
	the National Natural
	Science Foundation of China (11671366, 12071440 and 11801378).

\end{document}